\documentclass[11pt]{article}

\newcommand{\documentdate}{12 I 2021}

\usepackage{pgfplots}
\pgfplotsset{compat=1.13}
\usetikzlibrary{spy}

\usepackage{a4,latexsym,amssymb,varioref,graphicx,color,rotating,amsmath,mystyle}
\usepackage[a4paper]{geometry}

\title{Exploiting problem structure in derivative free optimization}

\author{M. Porcelli\thanks{%
Universit\`a di Bologna, Dipartimento di Matematica, 
Piazza di Porta S. Donato,5, 40126, Bologna, Italy.
Email: margherita.porcelli@unibo.it}\ \thanks{%
Institute of Information Science and
            Technologies ``A. Faedo'', ISTI-CNR, Pisa, Italy.}\,
~and  Ph. L. Toint\thanks{%
Namur Center for Complex Systems (NAXYS),
University of Namur,
61, rue de Bruxelles, B-5000 Namur, Belgium.
Email: philippe.toint@unamur.be}}
\date{\documentdate}

\begin{document}

\maketitle

\begin{abstract}
A structured version of derivative-free random pattern search optimization
algorithms is introduced which is able to exploit coordinate partially
separable structure (typically associated with sparsity) often present in
unconstrained and bound-constrained optimization problems.  This technique
improves performance by orders of magnitude and makes it possible to solve
large problems that otherwise are totally intractable by other derivative-free
methods. A library of interpolation-based modelling tools is also described,
which can be associated to the structured or unstructured versions of the
initial pattern search algorithm. The use of the library further enhances
performance, especially when associated with structure. The significant gains
in performance associated with these two techniques are illustrated using a
new freely-available release of the BFO (Brute Force Optimizer) package 
firstly introduced in \cite{PorcToin17}, which incorporates them. An interesting
conclusion of the numerical results presented is that providing global structural
information on a problem can result in significantly less evaluations of the
objective function than attempting to building local Taylor-like models.
\end{abstract}
{\bf Keywords:} derivative-free optimization, direct-search methods,
structured problems, interpolation models.\\

\noindent
{\bf Mathematics Subject Classification:} 65K05, 90C56, 90C90.

\numsection{Introduction}

Derivative-free methods have enjoyed a continued popularity and attention from
the very early days of numerical optimization
\cite{BoxWils51,HookJeev61,NeldMead65,Rose60} to nowadays \cite{AudeHare17,ConnScheVice08,LarsMeniWild19}.
 As their name indicates, such methods are aimed at solving nonlinear optimization 
problems without using derivatives of the objective function, which can be either too costly to
compute, or simply unavailable (as is often the case in simulation-based
applications). Among all possible methods, we focus here on two broad classes: model-based techniques
and pattern-search algorithms.  In the first class, an explicit local model of
the objective function is built, which is then used for finding any better
approximation of the sought minimizer. Early proposals in this direction
focused on unconstrained problems and include
\cite{ConnToin96,GripRina14,LuciScia02,Powe94a,Powe94b,Powe96,Powe98b,Winf73}, and were later
analyzed and/or extended to the constrained case (see
\cite{CartFialMartRobe19,ConnScheToin97,ConnScheVice03,ColsToin05,ConnScheVice09,CustRochVice10,Powe09,GratToinTroe10,ScheToin10,SampToin15,RobeCart19} among many others).
The second class, also initially for unconstrained problems
\cite{AbraAudeChriWals09,AudeDenn01,BoxWils51,DennTorc91,HookJeev61,NeldMead65,NeumFendSchlLeit11,Rose60}, was later analyzed in
\cite{CoopPric00,CoopPric01,CoopPric02,DolaLewiTorc03,GratRoyeViceZhan15,LagaPoonWrig12,Torc91,Torc97}, and also extended to more general contexts (see
\cite{AudeLeDiTrib19,LeDi11,LewiTorc00,PorcToin17,PricToin06,LiuzLuciRina20} for example).

\comment{
Arnold Neumaier, Hannes Fendl, Harald Schilly & Thomas Leitner 
\cite{CusRocVic10,CusVic07,GripRina15,
PhamWila11, LagaPoonWrig12,
CartGoulToin12a, ColsToin01a, ColsToin01b,  ConnScheToin97,
ConnScheVice08d, GratToinTshi11, LuciSciaTsen02, MaraNoce02,
Oeuv05, OeuvBier08, Powe98c,SampToin15
Powe65,ConnScheVice08,ConnScheVice08b,ConnScheVice08d}
}

Most of the contributions on derivative-free methods focus on small-scale
unstructured problems.  Indeed, building large-scale multivariate models is
very expensive as it requires a number of data points (and costly function
evaluations) which quickly grows with dimension: for the widely used quadratic
polynomial models, this number grows like the square of the problem's
size. As a consequence, model-based approaches are in general unrealistic
for moderate or large-scale applications. Similarly, the cost of the
sampling strategies at the heart of pattern-search techniques also explodes
with problem dimension and in practice restrict the application of these
algorithms to optimization in very few variables (a few tens).

The situation fortunately improves if, instead of pure black-box
optimization, one now considers the ``gray box'' case where one is allowed to
exploit some underlying problem structure (while still avoiding the use of any
derivative), see e.g. for an application \cite{MariMoriPorc18}.
Known sparsity in the objective function's Hessian was considered
in \cite{ColsToin01a,ColsToin01b} for methods based on quadratic models, where
it was shown that the growth in function evaluations with size is, for a large
class of sparse problems, essentially linear rather than quadratic.  However,
knowledge of the Hessian sparsity pattern is rarely directly obtained in
practice without the knowledge of derivatives. What is much more common,
for instance in discretized problems, is 
that optimization is performed on an application involving interconnected and
possibly overlapping subsystems. While the detailed analytic expression of the
subsystems' models are often unavailable, it is not unusual for their
connectivity pattern to be known. Such problems then often fall in the class
of \emph{coordinate partially separable} (or CPS) problems, where one 
attempts to solve
\beqn{pbps}
\min_{ x \in \Re^n} f(x) = \sum_{i=1}^q f_i(x)
\eeqn
where each $f_i: \Re^{n_i} \rightarrow \Re$ is a possibly nonsmooth (or even non
continuous) function depending on a subset $\{x_i\}_{i \in \calX_i}$ of the vector
of variables $x$, for some index sets $\calX_i \subseteq \{1, \dots, n\}$
($i = 1, \dots, q$). The cardinality of these $q$ subsets is denoted by
$n_i$ and is typically much smaller than the dimension $n$ of $x$.
In practice it is indeed common to see values of $n_i$ of
the order of ten or less, even for problems involving thousands of variables
or more, see e.g. problems described in Section~\ref{numerics_s}. 
Knowing the subsystems' connectivity pattern then typically amounts
to knowing the sets $\calX_i$ and being able to evaluate the ``element functions''
$f_i$ individually.  While we have introduced the concept for the unconstrained
formulation \req{pbps}, this is by no means restrictive, as constraints may be
added without affecting the structure.  In the rest of this paper, we will
assume that the problem's variables are subject to simple bounds, that is
\beqn{bounds}
\ell \leq x \leq u,
\eeqn
where $\ell$ and $u$ are vectors in $\Re^n$ such that $\ell \leq u$, all
inequalities being understood componentwise.

Optimizing coordinate partially separable
problems\footnote{Coordinate partially separable problems are a subclass of
  the more general partially separable problems, where the objective function
  takes the form $f(x) = \sum_{i=1}^q f_i(U_ix)$, where the $U_i$ are $n\times n$ low-rank
  matrices.} with model-based algorithms was considered in
\cite{ColsToin05}, where it was shown that the requirements in terms of
function evaluations now depend on the maximum value of the $n_i$, which is
often independent of the true problem's dimension\footnote{ Coordinate separability has been extensively considered 
also in the context of derivative-based algorithms, see e.g. \cite{RichTaka16,Wrig15}, but, to our knowledge,
in this setting the problem structure is assumed to be unknown and therefore is not exploited in the algorithm development.}.
The use of the same problem structure was also initiated in \cite{PricToin06}, where a first
pattern-search algorithm was described which could solve respectable size problems with $n$ up to $5625$.

The purpose of the present paper is to build on these contributions and to propose
(in Section~\ref{pattern-s}) a new derivative-free algorithm which is able to
exploit the problem structure at pattern-search level (known as a ``poll
step'') as well as a library of partially separable multivariate polynomials
to be combined with the structured poll step in a user-controlled ``search step'' (in
Section~\ref{model-s}).  The dramatic impact of these features on numerical performance
is then illustrated using a new version of the BFO (Brute Force Optimizer) package \cite{PorcToin17} incorporating them
(in Section~\ref{numerics_s}).  In particular, \emph{it is shown that the use of
the \emph{global} structural CPS information is often more efficient than
that of \emph{local} approximations of the objective function's Taylor expansion.}  It
also shown that \emph{a combined approach may also be advantageous,} at least
for problems of moderate size. Some conclusions and perspectives are finally
presented in Section~\ref{concl_s}.

\subsubsection*{Notation} The $j$-th component of a vector $x$ is denoted by
$x_j$, while $x_\calI$ denotes the subvector defined by considering the
components of $x$ indexed by the index set $\calI$.

\numsection{Exploiting coordinate partial separability in the poll step}
\label{pattern-s}

For the purpose of introducing our structure-exploiting techniques, we consider
 a generic pattern-search method consisting of 
a succession of iterations, each containing a search step (discussed in
Section~\ref{model-s}), a poll step, in which new (hopefully better) functions
values are generated by taking steps along randomly generated orthogonal sets
of search directions, and an update of the stepsizes as the iterations proceed.
A description of such a simple pattern search for unconstrained problems is given in Algorithm~\ref{simple-BFO}.
Comments on the constrained case are postponed to the end of this section.
\algo{simple-BFO}{Outline of a simple pattern search algorithm}{
  
{\bf Initialization}: The initial iterate $x$, the initial stepsize $\alpha$
and a convergence threshold $\epsilon>0$ are given. The parameters $\gamma \ge 1, \beta, \eta \in (0,1)$ are given. Define a set of orthonormal polling directions $\{d^i\}_{i=1}^{n}$. \\
  
 {\bf Until convergence}
  \begin{description}
  \item[Search step: ] Ask the user to provide a new (potentially improved) approximate
    minimizer of $f$, typically using problem specific modelling techniques;
  \item[Poll step: ] For $j = 1,\ldots,n$ or until ``sufficient decrease'' in $f$ is obtained 
    \begin{itemize}
      \item define a step of the form $\alpha d^j$
      \item evaluate $f(x+\alpha d^j)$ and, if necessary, $f(x-\alpha d^j)$,
    \end{itemize}
    If $f(x+\alpha d^j) < f(x)$ or $f(x-\alpha d^j) < f(x)$, replace $x$ by $x\pm\alpha d^j$,
    depending on which step gave decrease;
    
   Terminate the poll step if sufficient decrease is obtained, i.e.\ if
   \beqn{sufficient-decrease}
   f(x\pm\alpha d^j) - f(x) < \eta \alpha^2.
   \eeqn
  \item[Termination step: ]   \
  \begin{itemize}
   \item If $f$ has decreased in the course of the poll step, increase the
     stepsize $\alpha$ by setting $\alpha \leftarrow \gamma \alpha$,
     generate a new random set of orthonormal polling
     directions $\{d^i\}_{i=1}^{n}$ (successful iteration) and start a new iteration; 
   \item otherwise check for convergence (unsuccessful iteration):  if
     \begin{equation}\label{noconv}
     \alpha > \epsilon, 
     \end{equation}
     decrease $\alpha$ by setting $\alpha \leftarrow \beta \alpha$,
     generate a new random set of orthonormal polling directions $\{d^i\}_{i=1}^{n}$ and start a new iteration; otherwise declare convergence.
  \end{itemize}
  \end{description}
}

\clearpage

A few comments are useful at this stage. We first note that when sufficient
decrease is obtained for one of the polling direction, as tested in
\req{sufficient-decrease}, the algorithm stops using the current set of
directions and directly  updates the stepsize before starting a new iteration.


When \req{noconv} fails, that is when
\begin{equation}\label{conv}
 \alpha \leq \epsilon, 
\end{equation}
it would be acceptable to terminate the optimization entirely.  However, it is
useful to continue the effort for finding a better point by applying a
user-defined number of additional poll steps, each using a new random
orthonormal basis $\{d^i\}_{i=1}^{n}$. If condition \req{conv} is met at every
such step, final convergence is then declared.

We assume, as indicated by \req{pbps}, that all
variables are continuous
and that the polling
directions are given by the canonical coordinate vectors $\{e^i\}_{i=1}^n$.

Two key observations allow to set up a strategy to exploit the partial
separable structure in \req{pbps} that yields a new improved algorithm
inspired by \cite{PricToin06}. The first is that if a step along $e^j$ is
made, then only a subset of the element functions $f_i$ will be affected and
thereby need recomputation: only the $f_i$ such that $j$ appears in $\calX_i$
are concerned. Secondly, if variable $x_k$ does not occur in any of the
$f_i$'s involving variable $x_j$, a step of the form $\alpha_j e^j+\alpha_ke^k$
($\alpha_j$ and $\alpha_k$ are stepsizes) can be computed
involving completely disjoint sets of element functions: those involving
variable $x_j$ and those involving variable $x_k$. Crucially, the cost of this
combined step is potentially much \emph{less} than an unstructured step
along direction $e^j$ where the complete $f$ would be evaluated.

\vspace{\baselineskip}
\noindent{\bf Example.\ }\emph{ As an example, consider the partial separable function 
\begin{equation}\label{example}
f(x_1, x_2, x_3, x_4,x_5) = f_1(x_1, x_2) +
f_2(x_2,x_3) + f_3(x_1,x_2,x_4,x_5)+f_4(x_4,x_5) + f_5(x_4,x_5). 
\end{equation}
Then, a step along $e_1$, say, only involves the element functions $f_1$ and
$f_3$ while, since variable $x_3$ does not occur in $f_1$ and $f_3$, steps
along $\pm e_3$ can be made without affecting computations along $e_1$. We say
that $x_1$ and $x_3$ belong to disjoint sets of variables. Remarkably, the
evaluation at points of the form $\alpha_1 e_1+\alpha_3 e_3$ is less costly than
a full function evaluation since it does not involve the evaluation of $f_4$ and $f_5$.}
\vspace{\baselineskip}

Exploiting these observations in a systematic way is the key to the new
structured poll step.  More specifically, once the user provides the sets
$\calX_i$ ($i=1,\ldots,q$) indicating the subset of $\ii{n}$ defined by
the components of $x$ which appear in the domain of the element function
$f_i$($i=1,\ldots,q$), then this
structure is automatically first pre-processed as described in Algorithm \ref{preprocess}.

\algo{preprocess}{Subspace structure detection}{
\begin{enumerate}
\item The structure  $\{\calX_i\}_{i=1}^q$  is first inverted in that sets
  $\{\calE_j\}_{j=1}^n$ are built such that
  \[
  i \in \calE_j \tim{ if and only if } j \in \calX_i.
  \]
  The sets $\calE_j$ contains the indices of all element functions
  $f_i$ involving variable $x_j$.
\item Using lexicographic sorting, sets of variables' indeces
  $\{\calI_k\}_{k=1}^r$ corresponding to identical lists of elements $\calE_j$
  are constructed, that is 
  \[
  j_1 \in \calI_k \tim{ and } j_2 \in \calI_k \tim{ if and only if }
  \calE_{j_1} = \calE_{j_2} \eqdef \calY_k.
  \]
  The indeces of the element functions in each $\calY_k$ are thus
  indistinguishable as far as their dependence of the problem's variables
  is concerned.
\item \emph{Independent} collections $\{\calC_h\}_{h=1}^t$ of $\calI_k$ are
  then constructed by applying a greedy algorithm, independence being
  understood as the property that
  \[
  \calI_{k_1} \subseteq \calC_h \tim{ and } \calI_{k_2} \subseteq \calC_h
  \tim{ if and only if } \calY_{k_1} \cap \calY_{k_2} = \emptyset.
  \]
\end{enumerate}
}

As in \cite{PricToin06}, it can be verified that, for each $h \in \ii{t}$,
\beqn{Mh}
\calM_h \eqdef \bigcup_{ k \mid \calI_k \subseteq \calC_h} \calY_k 
\subseteq \ii{q}.
\eeqn
Thus steps in variables belonging to different sets $\calI_k$ in $\calC_h$ can be
computed independently (and in parallel).  Moreover, this process only involves a
subset $\calM_h$ of the set of all element functions.
Overall, the following collections are generated as a result of the procedure
described in Algorithm \ref{preprocess}:  $\{\calX_i\}_{i=1}^q$ and
$\{\calI_k\}_{k=1}^r$ containing variable indeces,  $\{\calE_j\}_{j=1}^n$ and
$\{\calM_h\}_{h=1}^t$ containing element function indeces and
$\{\calC_h\}_{h=1}^t$ containing indeces of sets of variables.

\vspace{\baselineskip}
\noindent{\bf Example.\ }\emph{
Returning to the the objective function described in \req{example}, we have that
$q = 5, n = 5$ and 
$$
\calX_1 = \{1,2\},\ \calX_2 = \{2,3\}, \ \calX_3 = \{1,2,4,5\},\ \calX_4 = \{4,5\},\ \calX_5 = \{4,5\}
$$
and the corresponding sets of element functions' indeces defined in Step 1 
of Algorithm \ref{preprocess} are given by
$$
\calE_1 = \{1,3\},\ \calE_2 = \{1,2,3\}, \ \calE_3 = \{2\},\ \calE_4 = \{3,4,5\},\ \calE_5 = \{3,4,5\}.
$$
Step 2 aims at finding repeated lists of element functions, i.e. $\calE_4 = 
\calE_5$ in this example, which are eliminated in the new sets  of variables'
and element functions' indeces. Therefore proceeding by lexicographic sorting,  we get $r=4$ sets
$$
\calI_1 = \{1\},\ \calI_2 = \{2\}, \ \calI_3 = \{3\},\ \calI_4 = \{4,5\},
$$
and
$$
\calY_1 = \{1,3\},\ \calY_2 = \{1,2,3\}, \ \calY_3 = \{2\},\ \calY_4 = \{3,4,5\}.
$$
We then define the subspaces of independent variables using a simple greedy
algorithm as follows. Initializing $\calC_1 = \{\calI_1\}$, we start checking
if $\calY_1 \cap \calY_2 = \emptyset$. Since,  $\calY_1 \cap \calY_2 \neq
\emptyset$ in our example, we  now check if $\calY_1 \cap \calY_3 =
\emptyset$. Since this is true, $\calI_3$ is added to $\calC_1$. Then we
check condition $(\calY_1\cup \calY_3) \cap \calY_4 = \emptyset$.
Since this condition fails, the definition of the set $\calC_1$ is completed
and the next set $\calC_2$ is initialized as $\calC_2= \{\calI_2\}$. Then one
checks the intersection of $\calY_2$ with $\calY_4$ ($\calI_3$ has already
been assigned). Since the intersection is not empty, we finally obtain $t = 3$
collections of sets of independent variables
\begin{equation}
  \label{groups}
\calC_1 = \{\calI_1, \calI_3\},\ \calC_2 = \{\calI_2\}, \ \calC_3 = \{\calI_4\}.
\end{equation}
The corresponding sets of  element functions' indeces are defined by
$$
\calM_1 = \calY_1 \cup \calY_3 = \{1,2,3\},\ \calM_2 = \calY_2 = \{1,2,3\}, \ \calM_3 = \calY_4 =\{3,4,5\},
$$
In this examples, steps for variables in each collection $\calC_i$, $i =
1,\dots, 3$ involves $|\calM_i| = 3$ element function evaluations. 
}
\vspace{\baselineskip}

Once the $\{\calC_h\}_{h=1}^t$ and $\{\calM_h\}_{h=1}^t$ are known, we may
define the subspaces
\begin{equation}\label{sk}
 \calS_k \eqdef \spanset{e^j \mid j \in \calI_k},
\end{equation}
$x_{\calS_k}$ the projection of the current iterate $x$ onto $\calS_k$ and 
the ``inactive'' index sets for the $h$-th collection as
\[
\calN_h = \ii{n} \setminus \bigcup_{k \mid \calI_k \subseteq \calC_h} \calI_k  
\ms\ms\ms \tim{for} h \in \ii{t}.
\]
The structured poll step then consists in performing poll steps in the
subspaces $\calS_k$  generated in each collection $\calC_h$ of sets of
independent variables and checking the sufficient decrease condition
\req{sufficient-decrease} in $f$ by only evaluating elements functions
in the corresponding set $\calM_h$ (for $k = 1, \dots, r$ and $h = 1, \dots, t$).

Our updated algorithm using a structured poll step is then given by Algorithm~\ref{BFOps}.

%

\algo{BFOps}{Pattern-search algorithm for partial separable problems}{
  
  {\bf Initialization}: The initial $x$,
  the initial stepsize $\alpha$, a convergence threshold $\epsilon>0$ 
  and the sets  $\calX_i \subseteq \{1, \dots, 
  n\}$ ($i = 1, \dots, q$) specifying the connectivity pattern   for problem
  \req{pbps} are given.
   The parameters $\gamma, \iota \ge 1, \beta, \eta \in (0,1)$ and $n_2 < n$ are given.  \\
  
  {\bf Structure analysis}: Compute a collection of sets of independent
  variables $\{\calC_h\}_{h=1}^t$ using Algorithm~\ref{preprocess}.
  Initialize the stepsizes $\alpha_k = \alpha$ for each set $\calI_k, k = 1, \dots, r$ and 
  define a set of orthonormal polling directions in $\calS_k, k = 1, \dots, r$ defined in \req{sk}.
  \\
  
 {\bf Until convergence}
  \begin{description}
  \item[1. Search step: ] Ask the user to provide a new (potentially improved) approximate
    minimizer of $f$, typically using problem specific modelling techniques;
  \item[2. Structured poll step: ] For $h = 1,\ldots, t$ or until ``sufficient decrease'' in $f$ is obtained
  \begin{description}
  \item[2.1 Poll step in $\calC_h$:]  For each $k$ such that $\calI_k \subseteq \calC_h$
    perform a standard poll step (using random orthogonal
    directions in $\calS_k$  and stepsize $\alpha_k$) starting from $x_{\calS_k}$ on the restricted
    function
    \[
    f_{\calY_k} = \sum_{i \in \calY_k} f_i,
    \]
    producing a potentially improved $x^+_{\calS_k}$ and decrease
    $f_{\calY_k}(x_{\calS_k}) - f_{\calY_k}(x^+_{\calS_k})\geq 0$;
    
    If sufficient decrease in $f_{\calY_k}$ is obtained, i.e. if
    $f_{\calY_k}(x^+_{\calS_k}) - f(x) < \eta \alpha_k^2$, increase the stepsize
    $\alpha_k$ by setting $\alpha_k \leftarrow \gamma \alpha_k$; otherwise
    decrease $\alpha_k$ by setting $\alpha_k \leftarrow \beta^{\iota} 
    \alpha_k$. Generate a new set of orthonormal polling directions in
    $\calS_k$. 
    
  \item[2.2 Iterate and function decrease update: ]  Define a new iterate $x^+$ by
    \[
      x^+_{\calI_k} =  x^+_{\calS_k} \tim{ for each } \calI_k \subseteq \calC_h
    \]
    and
    \[
      x^+_{\calN_h}  = x_{\calN_h}  \tim{otherwise,}
    \]
    and the corresponding objective function decrease by
    \[
    f(x^+) - f(x) = \sum_{k \mid \calI_k \subseteq \calC_h}
                    \Big[ f_{\calY_k}(x_{\calS_k}) - f_{\calY_k}(x^+_{\calS_k})\Big].
    \]
    If sufficient decrease in the objective function, i.e.\ if
    $f(x^+) - f(x) < \eta \alpha_k^2$,
    has been obtained, replace $x$ by $x^+$ and terminate the structured poll step.
  \end{description}
  \item[3. Termination step: ] Update the value of the ``global'' stepsize
    $\alpha = \min_k \alpha_k$. If sufficient decrease in $f$ was not obtained
    in the structured poll step and $\alpha$ is not below the prescribed
    accuracy $\epsilon$, go to Step 1. Otherwise enter the second pass (Step 4). 
    \item[4. Second pass: ] Generate a set of orthonormal polling directions
      $\{d^i\}_{i=1}^{n_2}$ in the full space and compute a complete poll loop
      along these directions with stepsize $\alpha$.  
    If sufficient reduction has been obtained, go to Step 1; else, declare convergence.
  \end{description}
}

\clearpage

As is standard in a poll step, an objective
function decrease may not be obtained in the structured poll step with the
current choice of stepsize, in which case
$x^+_{\calS_k}=x^{}_{\calS_k}$ and
$f_{\calY_k}(x_{\calS_k})=f_{\calY_k}(x^+_{\calS_k})$.

Once again we stress that using the successive collections $\{\calC_h\}_{h=1}^t$
is not mandatory (even if it clarifies the overall evaluation cost) and that
poll steps can be performed on the subspaces $\{\calS_k\}$ using the
restricted functions $\{f_{\calY_k}\}$ completely independently and in parallel.

Compared with the method described in \cite{PricToin06}, the algorithm
described in Algorithm~\ref{BFOps} is both simpler and more efficient, as it performs complete
poll steps using random orthogonal directions in each of the subspaces
$\calS_k$ (instead of only retaining the best increment for steps along a
fixed positive basis). Moreover, the mechanism used to adapt the
stepsizes is more elaborate than the somewhat adhoc
technique described in \cite[Section~2.3]{PricToin06}. It was also
observed in practice that the stepsize reduction can be slightly faster
than in the standard unstructured case. This faster reduction is
translated in the algorithm by a stepsize shrinking factor $\beta\in(0,1)$ (used in 
the unstructured case) which is raised to some power $\iota \geq 1$ (see Step 2.1 of Algorithm \ref{BFOps}). 

While the structured poll step is extremely efficient (as will be seen in the
example below and in Section~\ref{numerics_s}), it still has a potential
drawback.  Because the random poll directions are constrained to remain in
each of the $\calS_k$, they are not random directions in the complete space
$\Re^n$, and convergence is obtained, following \cite{GratRoyeViceZhan15}, to
points where no further decrease can be obtained in each of these
subspaces. This unfortunately does \emph{not} imply that no further decrease
can be obtained for $f$. A second pass is therefore necessary for obtaining
this desirable property (Step 4 in Algorithm \ref{BFOps}).  This second pass
does not use structure and is therefore considerably less efficient.  However,
this is mitigated by the observation that most (and sometimes all) the
decrease in objective function value is obtained in the first pass, the second
pass often only playing the role of a (possibly mildly expensive) convergence
check. Moreover, if sufficient decrease is identified during the second pass,
a return is made to the structure-using mechanism of the first, in an attempt
to efficiently improve the decrease obtained. Fortunately, the overall
efficiency of the structured minimization remains, in all examples we have
seen, orders of magnitude better than that of the unstructured one.  Moreover,
it is not unusual for applications of derivative-free algorithms that the user
is above all interested in obtaining a significant decrease in the objective
function and not so much in extracting the last carat of decrease, let
alone in checking local optimality.  In this case, the second unstructured
pass may often be unnecessary, bringing the optimization cost further down.

In order to reduce the cost of the second pass, our implementation of
Algorithm~\ref{BFOps} allows the user to specify a small number $n_2$ of random 
directions (typically much lower than $n$) for the second pass.

\vskip 7pt
We now discuss the cost of applying the structured poll step of
Algorithm~\ref{BFOps} compared with that of using Step~2 of
Algorithm~\ref{simple-BFO}.  Consider an unsuccessful poll step first, that is
a poll step during which sufficient decrease is not obtained.  (Note that such
steps must occur as the stepsizes have to become sufficiently small for the
algorithm to terminate.)  Because of \req{Mh}, we see that, for a given $h$,
the cost of evaluating the element function $f_i$ for $i \in \calM_h$ once cannot exceed
that of evaluating $f$, which we denote by $c_f$.  As forward and backward moves are considered for
every polling direction (at an unsuccessful poll step), the total evaluation cost of the
unsuccessful structured poll step is at most than $2tc_f$. By comparison, the
cost of an unsuccessful unstructured poll step is equal to $2nc_f$.
Since it is very often the case that $t \ll n$ (see Table \ref{probs-t}),
\emph{the evaluation cost of the structured step is typically only a small 
fraction of that of the unstructured one}.
Because the poll step is terminated as soon as sufficient decrease is
obtained, the cost of successful structured and unstructured poll steps is slightly more
difficult to compare.  As discussed in \cite{GratRoyeViceZhan15}, the expected
number of polling directions considered in a single successful (unstructured)
poll step is small (typically 2 or 3). Two variants are however possible
for the structured step. In principle, it can be terminated as soon as
sufficient decrease is obtained on a given subspace $\calS_k$. In the implementation 
discussed in Section~\ref{numerics_s}, the loop on the subspaces associated to a collection
$\calC_h$ is always completed before sufficient decrease triggers poll-step
termination. This choice appears to be efficient and allows for parallel
execution of the subspace-restricted poll steps for different subspaces.  Its
evaluation complexity therefore depends on the number of subspace collections
$\calC_h$ examined (which is at most $t$).  If this number is also a very
small integer (as is often the case), the evaluation cost of the structured
poll step is very similar to that of the unstructured one.  However, the
\emph{objective function decrease obtained is typically much larger}, as it
corresponds to a number of unstructured successful poll steps given by the
total number of subspaces $\calS_k$ considered in the calculation. Of course,
further gains may be obtained in the structured case if the parallelism
between the subspaces $\calS_k$ is also exploited.\\

When bound constraints are present, the backward and forward steps within the
poll step are truncated to prevent computation of the objective function at
infeasible points, but the rest of the calculation is essentially unchanged.

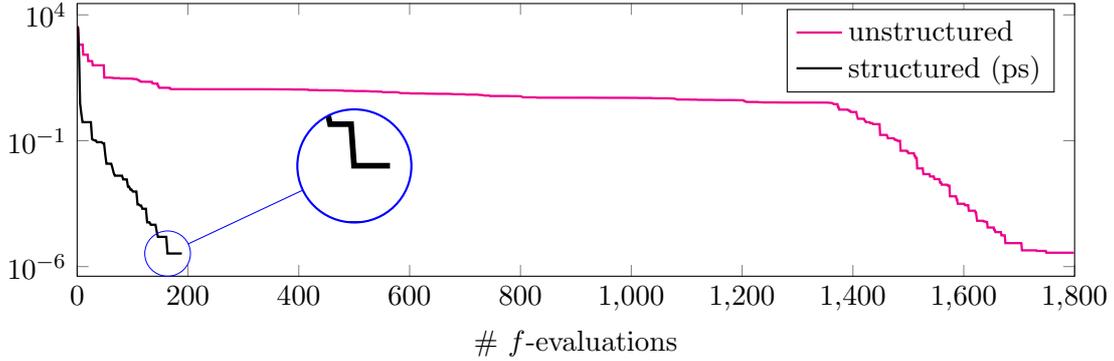
\begin{figure}[h]
\centering

\begin{tikzpicture}[spy using outlines=
	{circle, magnification=2.5, connect spies}]
\begin{semilogyaxis}[width=\linewidth, 
         xlabel = {\# $f$-evaluations},
         legend cell align=left,
         height = .25 \textheight,
         xmax = 1800,
         xmin = 0
         ]
 \addplot+[line width= 0.3mm,mark=none, magenta] table {./example_n10_nops.dat};
	    \addplot+[line width= 0.3mm,mark=none, black] table {./example_n10_ps.dat};  

  \coordinate (spypoint) at (axis cs:163,3.3E-6);
  \coordinate (magnifyglass) at (axis cs:500,0.01);
	\legend{unstructured, structured (ps)};
	\end{semilogyaxis}

\spy [blue, size=1.5cm] on (spypoint)
   in node[fill=white] at (magnifyglass);
\end{tikzpicture}

\caption{\label{fig:example}The evolution of $f$ as a function of the number
  of (complete) evaluations for the objective defined in
  \req{ex_f1}-\req{ex_f5}. Zoom on the function evaluations involved in the
  second pass when the partially separable structure is exploited (black curve).}
\end{figure}

\vspace{\baselineskip}
\noindent{\bf Example.}\emph{ We now give a first taste of the effectiveness of the structured poll step.
Consider an unconstrained problem with objective function structured as in
\req{example} and whose element functions are given by
\begin{eqnarray}
 f_1(x) & = & \sqrt{x_1^2+x_2^2}, \label{ex_f1}\\ 
 f_2(x) & = &  (\sin{x_1}-23 x_2 x_1)^2, \label{ex_f2}\\
 f_3(x) & = & (x_1^3-56 x_2 x_3)^2-x_4, \label{ex_f3}\\
 f_4(x) & = & \max\{|x_4|,|x_5|\}, \label{ex_f4}\\
 f_5(x) & = &  \sqrt{x_4^2+x_5^2}. \label{ex_f5}
\end{eqnarray}
We observe that the full objective $f$ is bounded below but the element $f_3$
is unbounded in variable $x_4$. However, the collection of independent sets in \req{groups} ensures 
that the restricted function $ f_{\calY_4}(x_{\calS_4}) $ is still bounded
below when computing poll steps in the subspace $S_4$ (the only one involving
$x_4$). In fact, the elements functions involved include $f_3$ together with
$f_4$ and $f_5$.
We then apply both the unstructured and structured variants of the method described above
with a tolerance $\epsilon$ on the stepsize equal to $10^{-5}$. Figure \ref{fig:example}
shows the two resulting convergence histories in terms of the number of
evaluations of the complete $f$. The advantage of using the structured poll step is striking.
}

\emph{The effect of the second pass is visible in the magnified part of the black curve (blue circle).
The first flat segment indicates that the second pass is entered, but
sufficient reduction is soon detected along a direction in the full space and the
structured poll step is performed again (vertical segment). Then, the second
pass is entered again (second flat segment) and convergence is finally declared.
}
\vspace{\baselineskip}

\numsection{A callable library of structured models}\label{model-s}

After discussing how to exploit the coordinate partially separable structure \req{pbps}
in the poll step, we now briefly describe a multivariate interpolation
technique that permits its exploitation in the user-controlled search step (i.e. Step~1 of
Algorithms~\ref{simple-BFO} and \ref{BFOps}). As is common for such steps (see \cite{LeDi11},
for instance), the idea is to provide a surrogate model of the objective
function in the neighbourhood of the current iterate $x_{best}$, which is
built using information gathered in the course of the algorithm.  This model
can then be minimized (typically within some trust region) to provide an
improved guess of the minimizer.  We now describe a library named
BFOSS for computing such a step while exploiting structure. 
As in \cite{ColsToin05,SampToin15},
it considers a surrogate model whose structure mirrors that of \req{pbps} in
that it shares the same coordinate partially separable definition.  This is
achieved by constructing separate multivariate polynomial interpolation models
for each of the $f_i$ (restricted to their free variables). General
multivariate polynomial interpolation models follow the principles of
\cite{ConnScheToin97,ScheToin10,ConnScheVice09}, while their use in the
context of partially separable problems is inspired by
\cite{ColsToin05,SampToin15} and the necessary adaptations used here for
handling the bound constraints are similar to those discussed in
\cite{GratToinTroe10}.  While the general idea of structured models is not
new, our experience suggests that its practical performance does depend on a
number of more detailed decisions at the implementation level.  We therefore focus, in what follows, on the algorithmic
aspects which are specific to our approach and relevant to the numerical
comparison to be conducted in Section~\ref{numerics_s}. 
We first describe our strategy for the
unstructured case ($q=1$) and then specialize it to the partially separable
case ($q \geq 1$).
\comment{
  and let $c_{mesh}$ be the current 
mesh size.
Let $Y_k = \{y_i\}$ be $p$ points (including $x_{best}$) and 
a trust-region be defined as  $ \{x_{best} + s: s \in 
\Re^n \mbox{ and } \|s\|_{\infty} \le \Delta\}$ with radius 
$\Delta = \min \{\delta c_{mesh}, \min(u-l)/2)\}$ for a factor 
$\delta \ge 1$.
}

Let $Y = \{y^i\}_{i=1}^p$ be a sample set and let $m(x)$ denote a
polynomial of degree $d$ ($d = 1,2$)  interpolating $f$ at the points in $Y$,
that is satisfying the interpolation conditions 
\begin{equation}\label{modcond}
m(y^i) = f(y^i), \ \  i = 1, \dots, p.
\end{equation}
The polynomial $m(x)$ can be expressed as a linear combination
of elements of the natural basis $\phi$ for the space of polynomial
of degree at most $2$, i.e. the basis of
monomials, that we assumed to be ordered as follows 
$$\phi =\left \{1,x_1,x_2, \dots,x_n,\frac{1}{2}x_1^2, \dots, \frac{1}{2}x_n^2, 
x_1x_2, \dots, x_{n-1}x_n, x_1x_3, \dots, x_{n-2}x_n, \dots, x_1x_n\right \}.$$
%

Finding a model that satisfies conditions \req{modcond} is
equivalent to solving a potentially underdetermined linear system of the form
\begin{equation}\label{polsys}
 M_{\phi} z = f(Y)
\end{equation}
with coefficient matrix $M_{\phi}$ of dimension $p \times   \bar p,
z \in \Re^{\bar p}$ and $(f(Y))_i=f(y^i), i = 1,\dots,p,$ where $p \le \bar p$
and $\bar p$ is the chosen number of elements from the basis $\phi$.

If the coefficient matrix $M_{\phi}$ is square and nonsingular, then the model $m$ is unique. 
If $p <  \bar p$ the linear system is underdetermined,  
the resulting interpolating polynomials may not exist or may no
longer be unique, and different approaches to construct
the model $m(x)$ are possible \cite{ConnScheVice09}.
BFOSS allows the user to choose between two possibilities: a model can be
built by taking the minimum 2-norm solution of \req{polsys}
\cite{CustRochVice10} or, alternatively, a sub-basis $\tilde \phi$ of  $p$
elements is extracted from $\phi$ yielding the corresponding square matrix
$M_{\tilde \phi}$ (the components of $z$ corresponding to the removed columns
are set equal to zero).
 
%
%
%


As a result of the search step, a new tentative iterate $x_+$ is then
computed by minimizing the model $m$ in the intersection between the
trust-region and the feasible set, that is
\beqn{trbox}
  \{x_{best} + s:\ \max(l-x_{best}, -\Delta)\le s \le
\min(u-x_{best}, \Delta) \}
\eeqn
where $\Delta>0$ is the current trust-region radius, and $f(x_+)$ is evaluated.
Computing $x_+$ within the trust-region \req{trbox} is a well-understood question
and there are several algorithms available for the task (see \cite[Chapter~7]{ConnGoulToin00} for an overview).
BFOSS allows the user to choose between a bound-constrained variant of the 
Mor\'e-Sorensen trust-region algorithm and a projected gradient
trust-region algorithm to solve \req{trbox} \cite{MoreSore83,CartGoulToin09a,Toin81b,Stei83a}.
The trust-region radius $\Delta$ is updated on the basis of the ratio of
the achieved to predicted reduction $\rho = (f(x_+)-f(x_{best}))/(m(x_+)-m(x_{best}))$ 
as follows:
\[
\Delta = \left \{\begin{array}{ll}
\min \{\alpha_1 \Delta, \Delta_{max}\} & \mbox{ if } \rho> \eta_2 \\
\Delta                               & \mbox{ if } \rho \in [\eta_1,\eta_2] \\ 
\max \{ \epsilon_m, \alpha_2 \|s\|\}   & \mbox{ else } 
\end{array} \right .
\]
where $\alpha_1>1, \alpha_2 \in (0,1)$, $\Delta_{max}>0$ and $0<\eta_1<\eta_2<1$ are
trust-region parameters \cite[Chapter 17]{ConnGoulToin00} and $\epsilon_m$ denotes the machine precision.
When BFOSS is called the first time in Step~1 of
Algorithm~\ref{simple-BFO}, the value of $\Delta$ is set equal to the current step size.
At all subsequent calls, if $\Delta$ is too small, i.e. 
$\Delta < \Delta_{min}$ with $\Delta_{min}$ defined by the  
user, the radius $\Delta$ is restarted to the current step size.
If $\Delta$ is below $\Delta_{min}$ after this restart, BFOSS is terminated.
Finally, all the newly evaluated points
($x_+$ and, when relevant, the new interpolation points) are returned
to Algorithm~\ref{simple-BFO} with their associated function values.

It is well-known that the fact that the model $m$ is well-defined not only depends on
the number of points in $Y$, but that a further geometric condition (known as
{\em poisedness}) is also required (see \cite{ConnScheVice09} for details).
We now describe the specific strategy used in BFOSS to define
the interpolation set $Y$, which combines the use of models of various
degrees and types and ensures the necessary poisedness condition.

When using BFOSS, the user selects a value of $\bar p$ in the set $\{n+1,
2n+1, (n+1)(n+2)/2\}$ corresponding to a linear model, a quadratic one with a
diagonal Hessian and a full quadratic one.
A first tentative set $Y$ of $\bar p$ points is built around $x_{best}$ extracting 
from the recorded history (past points and function values), the closest points to $x_{best}$ in 
the 2-norm  and their associated function values (including $x_{best}/f(x_{best})$).
If the history is too short in that it only contains $p<\bar p$ points, BFOSS
uses all the available $p$ points according to one of the strategy described
above for the case $p<\bar p$.

Then matrix $M_{\phi}$ in \req{polsys}  is then built, 
together with its generalized inverse $M_{\phi}^\dagger$. 
The generalized inverse $M_{\phi}^\dagger$ is constructed from a 
truncated SVD of $M_{\phi}$, where  ``redundant'' singular values in a
standard SVD are zeroed. We say that a singular value is 
redundant if its value is lower than the maximum singular value scaled by a user-defined parameter
$k_{ill}$ that measure the maximum ill-conditioning allowed.
\footnote{By default $k_{ill}= \infty$ and $k_{ill}= 10^{12}$ for structured and unstructured 
problems, respectively. These values resulted from our numerical experiments.}
Finally, the pseudoinverse is computed from this regularized SVD. 
When redundant singular values are detected, the corresponding points in $Y$
are progressively removed, replaced by random points and $M_{\phi}^\dagger$
rebuilt. Importantly, this crucial phase does not require the computation of
additional function values. 

Once $M_{\phi}$ is considered sufficiently well-conditioned (within  $k_{ill}$)
the poisedness of the current set $Y$ is measured by computing the maximum 
absolute value of the Lagrange  polynomials in the neighbourhood.
Based on this measure, following \cite{SampToin16},
some points in $Y$ can be replaced by some ``far'' points available 
in the history but not used so far.
Then, $Y$ is possibly further improved by performing exchanges until the improvement
in poisedness becomes moderate enough. Finally, the obtained 
(reasonably conditioned)  $M_{\phi}^\dagger$ matrix is used to define the linear
combination of the Lagrange polynomials which interpolates function values at the
interpolation points.

\comment{
We note that by default, if the search step is used in BFO, all evaluated 
points are recorded and therefore they are used to form $Y$ after the few 
first iterations. On the other hand, if the user chooses not to save any 
information for memory reasons ($\ell =1$), each interpolation set is formed, used 
and then discarded.
}

The adaptation of this strategy to the partially separable 
case is straightforward: one simply applies the same technique to define
interpolation models $m_i$ for each of the element functions $f_i(x_i)
(i = 1, \dots , q)$.
\comment{
Therefore it is easy to define $q$ interpolation sets $Y_i(i = 1,\dots, q)$
of dimension $p_i$ around each $x_{best}$ restricted to the $i$th domain 
$\Re^{n_i}$ and corresponding interpolation models $m_i$ 
\begin{equation}\label{modeli}
 m_i(x_{best}+s) = c_i + g_i^T s + \frac{1}{2}s^T H_i s, \quad i = 1, \dots, q,
\end{equation}
where $c_i\in\Re, g_i\in\Re^{n_i}$ and $H_i\in \Re^{n_i \times n_i}, \
i = 1, \dots, q$.
}
Note that each model $m_i$ has at most $n_i$ variables and 
approximate $f_i$'s around the projection of $x_{best}$ onto the 
subspace $\Re^{n_i}$. 
The final trial point $x_+$ is then
computed by minimizing the global quadratic model 
$$
m(x_{best}+s) = \sum_{i=1}^q m_i(x_{best}+s), 
$$
in the box \req{trbox}.

\numsection{Numerical illustration}\label{numerics_s}

We now report numerical experiments comparing and combining the two techniques
described above.  The results were obtained by implementing the structured
poll-step described in Algorithm~\ref{BFOps} with the BFO (the Brute Force Optimizer) 
package \cite{PorcToin17} and exploiting the BFOSS library in combination with this upgraded version of BFO\footnote{The BFOSS
library is distinct from the BFO package itself --they come in different
files-- although it interacts with it.  But its use or even presence is not
necessary for running the main package.}. BFO is a random
pattern search algorithm proposed by the authors in \cite{PorcToin17}, whose
structure is identical to that of Algorithm~\ref{simple-BFO}.  Since a
full description of this package is somewhat involved and many of its features
irrelevant for our present discussion, we avoid restating it here in detail,
and refer the reader to \cite{PorcToin17} for an in-depth description.

The main objectives of the numerical
tests are
\begin{itemize}
\item to illustrate the impact of structure usage both at rather global level
  (via the use of the structured poll step) and a a local level (by building
  local interpolation/regression models using BFOSS),
\item to discuss the relative merits of the two techniques,
\item to investigate the effectiveness of their combination.
\end{itemize}

We first provide some numerical illustration of our claim that the structured
poll step is more efficient than the unstructured one.  For this purpose, we
compare the structured\footnote{In the BFO implementation of the new poll step,
$\iota$ is a newly introduced algorithmic parameter which, like all such
parameters \cite{PorcToin17}, can be (and has been) trained for improved
performance.} and unstructured version of BFO\footnote{It may be
recalled that the unstructured BFO was shown in \cite{PorcToin17} to be
quite competitive, in particular when compared with NOMAD \cite{LeDi11}.} on
a set of variable dimension test problems extracted from {\sf CUTEst}
\cite{GoulOrbaToin15b} and/or already used in \cite{PricToin06} for the most
part.

For each problem, we considered dimensions ranging from around 10 to 10000, whenever
solvable in reasonable time\footnote{Eight hours using Matlab
R2017b on a Intel(R) Xeon(R) CPU E3-1245 v5 @ 3.50GHz machine with 64 GB RAM.}. 
They are partitioned in four dimension-dependent test sets: small, smallish, medium and large.
Variables are all considered to be continuous although the 
structured poll step is implemented in BFO for handling integer and categorical variables as well.
The performance analysis  with these variables are out of the scope of this work and will be considered in the future.
Their main characteristics are detailed in Table~\ref{probs-t}, but we now 
give some additional information.
\begin{itemize}
\item The original {\tt BEALE} problem from {\sf CUTEst} only has two
  variables.  The variant {\tt BEALES} used here is obtained by juxtaposing $n/2$ copies
  of the original problem, resulting in a totally separable problem (hence the {\tt S}).
  It is useful as it exposes how well a method can exploit such an important structure.
\item The function {\tt NZF1} is derived from that published in
  \cite{PricToin06}.  Its precise formulation is detailed in Appendix \ref{app:nzf1}. It
  reduces to the version of \cite{PricToin06} for $n=13$.
\item The {\tt CONTACT} problem is a bound-constrained minimum surface problem
  involving nonlinear surface boundaries and an obstacle from below the
  surface.  It is described in Appendix \ref{app:contact}.
\item The {\sf CUTEst} minimum-surface problem {\tt NLMINSRF} only differs
  from {\tt LMINSURF} (also in {\sf CUTEst}) in that the surface boundary is
  nonlinear.
\item The {\tt BROWNAL6} problem is that presented in \cite{PricToin06} under the name
  'Brown almost linear'.
\item In the original {\tt MOREBV} problem, the distance from the starting
  point to the solution decreases with dimension, which makes it less
  interesting for large problems. We use here a version of the problem where
  the original starting point is multiplied by $\log_{10}(n)$ to compensate.
\item The {\tt JNLBRNG1}  is the
    journal bearing problems of MINPACK2 \cite{AverMore92}. 
\item Problems {\tt BROYDN3D}, {\tt CONTACT}, {\tt ENGVAL}, {\tt JNLBRNG1}, {\tt LMINSRF}, {\tt NLMINSRF}, {\tt MOREBV}
are discretized problems. They illustrate how frequently the CPS structure appears in the situation: 
one variable associated to a particular location in the problem
typically only depends on the variables associated with neighbouring locations, and not with all of them.
\end{itemize}
These problems are typical of a very large class of medium/large-scale applications, where, although
$q$ obviously depends on $n$, both $t$ and $\max_k|\calI_k|$, the maximal dimension of any subspace
$\calS_k$, do not.

\begin{table}
\begin{center}
\scriptsize
\begin{tabular}{lcccccccc} 
Problem    &\multicolumn{4}{c}{instance's dimensions ($n$)} & $q$ & $\max_k|\calE_k|$ & $t$ & $\max_k|\calI_k|$\\
           & small & smallish & medium & large       &                   &                   &      &          \\
\hline
           &                                         &                   &                   &      &          \\
ARWHEAD    & 10& 50, 100& 500, 1000& 5000, 10000     &  $n-1$            &         2         &   2  &       1  \\
BDARWHD    & 10& 50, 100& 500, 1000& 5000, 10000     &  $n-2$            &         3         &   3  &       1  \\
BDEXP      & 10& 50, 100& 500, 1000& 5000, 10000     &  $n-2$            &         3         &   3  &       1  \\
BDQRTIC    & 10& 50, 100& 500, 1000& 5000, 10000     &  $n-4$            &         5         &   5  &       1  \\
BEALES     & 10& 50, 100& 500, 1000& 5000, 10000     &  $n/2$            &         2         &   1  &       2  \\
BROYDN3D   & 10& 50, 100& 500, 1000& 5000, 10000     &  $n-1$            &         3         &   3  &       1  \\
BROWNAL6   & 10& 50, 102& 502, 1002& 5002, 10002     &  $(n-2)/4$        &         6         &   2  &       4  \\
CONTACT    & 15& 64, 144& 400, ~\,900& 2500, ~\,4900 &  $(\sqrt{n}-1)^2$ &         4         &   4  &       1  \\
ENGVAL     & 10& 50, 100& 500, 1000& 5000, 10000     &  $n-1$            &         2         &   2  &       1  \\
DIXMAANA   & 15& 51, 102& 501, 1002& 5001, 10002     &  $n$              &         4         &   4  &       1  \\
DIXMAANI   & 15& 51, 102& 501, 1002& 5001, 10002     &  $n$              &         4         &   6  &       1  \\
FREUROTH   & 10& 50, 100& 500, 1000& 5000, 10000     &  $n-1$            &         2         &   2  &       1  \\
HELIX      & 11& 21, 101& 501, 1001& 5001, 10001     &  $(n-1)/2$        &         3         &   2  &       2  \\
JNLBRNG1   & 24& 64, 144& 400, ~\,900& 2500, ~\,4900 &  $4(n_t+1)^2$       &         3         &   3  &       1  \\
           &   &        &            &               &  $(n=2(n_t+2)(n_t+1))$ &                   &      &          \\
LMINSURF   & 16& 64, 144& 400, ~\,900& 2500, ~\,4900 &  $(\sqrt{n}-1)^2$ &         4         &   4  &       1  \\
MOREBV     & 12& 52, 102& 502, 1002& 5002, 100002    &  $n$              &         3         &   3  &       1  \\
NLMINSRF   & 16& 64, 144& 400, ~\,900& 2500, ~\,4900 &  $(\sqrt{n}-1)^2$ &         4         &   4  &       1  \\
NZF1       & 13& 39, 130& 650, 1300& 6500, 13000     &  $(7n/13)-2$      &         6         &   4  &       2  \\
POWSING    & 20& 52, 100& 500, 1000& 5000, 10000     &  $n/4$            &         4         &   1  &       4  \\
ROSENBR    & 10& 50, 100& 500, 1000& 5000, 10000     &  $n/2$            &         2         &   1  &       2  \\
TRIDIA     & 10& 50, 100& 500, 1000& 5000, 10000     &  $n$              &         2         &   2  &       1  \\
WOODS      & 20& 40, 200& 400, 2000& 4000, 10000     &  $n/4$            &         4         &   1  &       4  \\
\end{tabular}  
\caption{Characteristics of the test problem instances}
\label{probs-t}
\end{center}
\end{table}

Our experiments were conducted using a new release of BFO, which contains both
the new structure-exploiting poll step and the BFOSS library. In order to
train the parameter $\iota$ defining the faster stepsize decrease in the
structured optimization pass, we selected, for each of the above problems, the
instance of third smallest dimension (mostly $n \approx 100$) and performed
training to improve the resulting data profile (see \cite{PorcToin17c}). The
experimental guess of $\iota = 1.25$, was only very marginally
improved\footnote{We used the BFO default training accuracy requirement
  $\epsilon_t=0.01$.} to $1.2550$. Note that other BFO algorithmic parameters
were set to their default values.  For small problems, the reported results
are an average of 30 independent runs, for smallish 10 runs for medium ones, 5
runs and a single run for the large ones.  All the results discussed in this
section are reported in Table~\ref{results-full-t} in Appendix \ref{app:res}.

\subsection{Exploiting structure in the poll step}

The results obtained by running the new structure-exploiting version of BFO
(with the trained $\iota$) are presented in Table~\ref{results-t}.  This table
reports the number of complete function evaluation to obtain an
approximate \footnote{We used the BFO default accuracy requirement $\epsilon =
0.0001$.  We also stress that the algorithmic parameters in BFOSS have been
trained for improved performance (as all other parameters of BFO) on the small
test-set, using the data profile performance measure \cite{PorcToin17c}.
The values of the trained parameters are set as default  in BFO and BFOSS.} 
minimizer for each of the instances of Table~\ref{probs-t}.  For each
instance, we give the number of required objective-function (full) evaluations
for structure exploiting (first) and standard (second, no structure
exploitation) versions of BFO \footnote{When partially separable structure is provided, 
the number of full function evaluation per BFO iteration is retrieved by summing the number of element function evaluations
and dividing  by the number of elements  (and rounding to integer). }.

\begin{table}
\begin{center}
\scriptsize
\begin{tabular}{l|ccccccc}
Problem      &    small       & \multicolumn{2}{c}{smallish}    &  \multicolumn{2}{c}{medium}    & \multicolumn{2}{c}{large}     \\
             &                &                &                &                &               &               &               \\
\hline
             &                &                &                &                &               &               &               \\
ARWHEAD      &    79/962      &     91/11859   &    97/36085    &   146/$\infty$ &  194/$\infty$ &  389/$\infty$ &  618/$\infty$ \\
BDARWHD      &   101/2152     &    70/$\infty$ &    71/$\infty$ &    81/$\infty$ &   76/$\infty$ &   77/$\infty$ &   77/$\infty$ \\
BDEXP        &  1661/21122    & 13218/$\infty$ & 34409/$\infty$ &$\infty$/$\infty$&$\infty$/$\infty$&$\infty$/$\infty$&$\infty$/$\infty$\\
BDQRTIC      &   290/2468     &   301/$\infty$ &   298/$\infty$ &   393/$\infty$ &  542/$\infty$ & 1150/$\infty$ & 1480/$\infty$ \\
BEALES       &   275/$\infty$ &   275/$\infty$ &   275/$\infty$ &   275/$\infty$ &  275/$\infty$ & 300/$\infty$ &  325/$\infty$ \\
BROYDN3D     &   308/1225     &   199/22244    &   273/74758    &   304/$\infty$ &  370/$\infty$ &  640/$\infty$ &  675/$\infty$ \\
BROWNAL6     &   15773/6171   &  4325/$\infty$ &  2788/$\infty$ &  2682/$\infty$ & 2847/$\infty$ & 3059/$\infty$ & 3118/$\infty$ \\
CONTACT      &   265/268      &   221/30452    &   442/$\infty$ &  1009/$\infty$ & 1761/$\infty$ & 3589/$\infty$ & 5952/$\infty$ \\
DIXMAANA     &   185/2358     &   189/17357    &   240/55748    &   226/$\infty$ &  375/$\infty$ &  310/$\infty$ &  756/$\infty$ \\
DIXMAANI     &   185/17702    &   186/$\infty$ &   184/$\infty$ &   230/$\infty$ &  265/$\infty$ &  709/$\infty$ &  798/$\infty$ \\
ENGVAL       &   143/1567     &   155/30400    &   157/$\infty$ &   159/$\infty$ &  159/$\infty$ &  159/$\infty$ &  159/$\infty$ \\
FREUROTH     &   257/83101    &   181/$\infty$ &   191/$\infty$ &   185/$\infty$ &  192/$\infty$ &  233/$\infty$ &  317/$\infty$ \\
HELIX        &   131/9931     &   152/$\infty$ &   167/$\infty$ &   298/$\infty$ &  388/$\infty$ &  883/$\infty$ & 1070/$\infty$ \\
JNLBNG1      &   101/734      &   301/27567    &   427/$\infty$ &   944/$\infty$ & 1393/$\infty$ & 1909/$\infty$ & 1799/$\infty$  \\
LMINSURF     &   306/318      &   461/23692    &  1065/$\infty$ &  3301/$\infty$ & 8506/$\infty$ & 29961/$\infty$ &$\infty$/$\infty$ \\ 
MOREBV       &  7010/7154     &    60/22732    &    47/1923     &    47/9001     &   47/18001    &   47/90001    &   47/$\infty$ \\
NLMINSRF     &   412/415      &   579/29409    &  1082/$\infty$ &  8853/$\infty$ & 3412/$\infty$ &31619/$\infty$ &$\infty$/$\infty$ \\ 
NZF1         &   177/1480     &   225/61772    &   625/$\infty$ &   684/$\infty$ &  667/$\infty$ &  946/$\infty$ & 2228/$\infty$ \\
POWSING      &   716/20605    &   824/$\infty$ &   849/$\infty$ &   988/$\infty$ & 1036/$\infty$ & 1164/$\infty$ & 1148/$\infty$ \\
ROSENBR      &   361/13241    &   361/$\infty$ &   384/$\infty$ &   436/$\infty$ &  461/$\infty$ &  636/$\infty$ & 736/$\infty$ \\
TRIDIA       &   440/3073     &   316/$\infty$ &   314/$\infty$ &   345/$\infty$ &  293/$\infty$ &  277/$\infty$ & 278/$\infty$ \\
WOODS        &  1609/$\infty$ &  1747/$\infty$ &  1924/$\infty$ &  2241/$\infty$ & 2927/$\infty$ & 4515/$\infty$ & 5002/$\infty$ \\
\end{tabular}
\caption{\label{results-t}The numbers of objective function evaluations required by
  the structured/unstructured versions of BFO for the problem instances of
  Table~\ref{probs-t}, as a function of increasing problem size ($\infty$
  meaning that more than 100000 evaluations were needed).} 
\end{center}
\end{table}

It clearly results from Table~\ref{results-t} that \emph{it is possible to
  solve large partially-separable problems without using
  derivatives at an acceptable cost in number of function evaluations.  Using
  the structure is crucial if the problem size exceeds ten or so}. In fact, a
significant fraction of the problems of that size can't be solved at all in a
reasonable number of evaluations if structure is neglected: direct
derivative-free methods like BFO proceed by sampling, and this technique badly
suffers from the curse of dimensionality.

A comparison of the above results with those of Table~2 in \cite{PricToin06}
also shows that the structured BFO significantly outperforms the simpler
algorithm presented in that reference.

We now illustrate our comments of Section~\ref{pattern-s} about the relative
efficiency of the two polling passes. Figures~\ref{br10-50-f} to
\ref{br1000-f} show, for the {\tt BROYDEN3D} problems in dimensions ten to one
thousand, the objective function decreases obtained by the new BFO using
structure (black line) and that obtained by the standard version of BFO
which ignores structure (magenta line).  These decreases are
expressed as a function of the number of (complete) objective-function
evaluations.

\comment{
\begin{figure}[htbp]
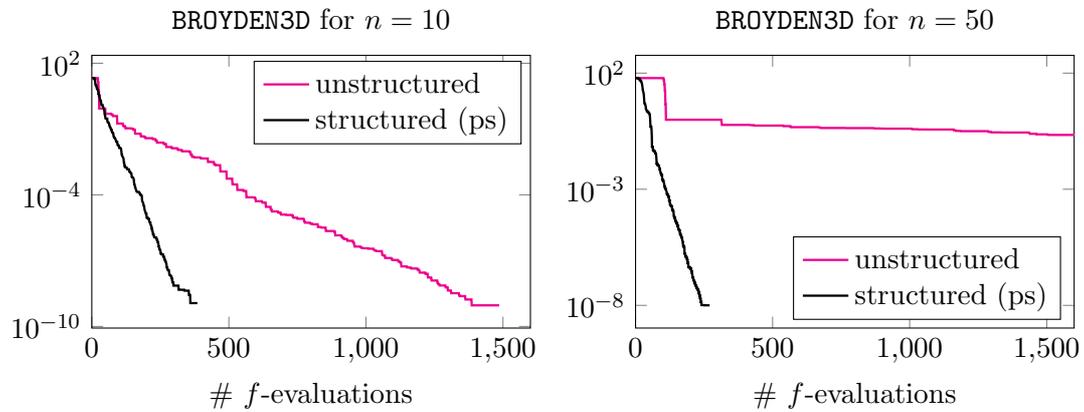

\centerline{
\includegraphics[height=5cm,width=7cm]{bd3d_n10.eps}
\includegraphics[height=5cm,width=7cm]{bd3d_n50.eps}
}
\caption{\label{br10-50-f}The evolution of $f$ as a function of the number of (complete)
  evaluations for {\tt BROYDEN3D} for $n=10$ and $n=50$.}
\end{figure}
\begin{figure}[htbp]
\centerline{
\includegraphics[height=5cm,width=7cm]{bd3d_n100.eps}
\includegraphics[height=5cm,width=7cm]{bd3d_n500.eps}
}
\caption{\label{br100-500-f}The evolution of $f$ as a function of the number of (complete)
  evaluations for {\tt BROYDEN3D} for $n=100$ and $n=500$.}
\end{figure}

\begin{figure}[htbp]
\centerline{
\includegraphics[height=5cm,width=7cm]{bd3d_n1000.eps}
}
\caption{\label{br1000-f}The evolution of $f$ as a function of the number of (complete)
  evaluations for {\sc BROYDEN3D} for $n=1000$.}
\end{figure}
}

\begin{figure} \centering
\begin{tikzpicture}
\begin{semilogyaxis}[width=.5 * \linewidth, 
         title = {{\tt BROYDEN3D} for $n=10$}, 
         xlabel = {\# $f$-evaluations},
         legend cell align=left,
         height = .25 \textheight,
         xmax = 1600,
         xmin = 0
         ]       
	    \addplot+[line width= 0.3mm,mark=none, magenta] table {./broy3d_n10_nops.dat};
	    \addplot+[line width= 0.3mm,mark=none, black] table {./broy3d_n10_ps.dat};    
	\legend{unstructured, structured (ps)};
	\end{semilogyaxis}
	\end{tikzpicture}
	\begin{tikzpicture}
        \begin{semilogyaxis}[width=.5 * \linewidth, 
         title = {{\tt BROYDEN3D} for $n=50$}, 
         legend pos=south east,
         xlabel = {\# $f$-evaluations},
         legend cell align=left,
         height = .25 \textheight,
         xmax = 1600,
         xmin = 0
         ]       
	    \addplot+[line width= 0.3mm,mark=none, magenta] table {./broy3d_n50_nops.dat};
	    \addplot+[line width= 0.3mm,mark=none, black] table {./broy3d_n50_ps.dat};	    
	\legend{unstructured, structured (ps)};
	\end{semilogyaxis}
	\end{tikzpicture}
	\caption{\label{br10-50-f}The evolution of $f$ as a function of the number of (complete)
  evaluations for {\tt BROYDEN3D} for $n=10$ and $n=50$}
	\end{figure}

\begin{figure} \centering
\begin{tikzpicture}
\begin{semilogyaxis}[width=.5 * \linewidth, 
         title = {{\tt BROYDEN3D} for $n=100$}, 
         legend pos=south east,
         xlabel = {\# $f$-evaluations},
         legend cell align=left,
         height = .25 \textheight,
         xmax = 1600,
         xmin = 0
         ]       
	    \addplot+[line width= 0.3mm,mark=none, magenta] table {./broy3d_n100_nops.dat};
	    \addplot+[line width= 0.3mm,mark=none, black] table {./broy3d_n100_ps.dat};    
	\legend{unstructured, structured (ps)};
	\end{semilogyaxis}
	\end{tikzpicture}
	\begin{tikzpicture}
        \begin{semilogyaxis}[width=.5 * \linewidth, 
         title = {{\tt BROYDEN3D} for $n=500$}, 
         legend pos=south east,
         xlabel = {\# $f$-evaluations},
         legend cell align=left,
         height = .25 \textheight,
         xmax = 1600,
         xmin = 0
         ]       
	    \addplot+[line width= 0.3mm,mark=none, magenta] table {./broy3d_n500_nops.dat};
	    \addplot+[line width= 0.3mm,mark=none, black] table {./broy3d_n500_ps.dat};	    
	\legend{unstructured, structured (ps)};
	\end{semilogyaxis}
	\end{tikzpicture}
	\caption{\label{br100-500-f}The evolution of $f$ as a function of the number of (complete)
  evaluations for {\tt BROYDEN3D} for $n=100$ and $n=500$}
	\end{figure}

\begin{figure} \centering
\begin{tikzpicture}
\begin{semilogyaxis}[width=.5 * \linewidth, 
         title = {{\tt BROYDEN3D} for $n=1000$}, 
         legend pos=south east,
         xlabel = {\# $f$-evaluations},
         legend cell align=left,
         height = .25 \textheight,
         xmax = 1600,
         xmin = 0
         ]       
	    \addplot+[line width= 0.3mm,mark=none, magenta] table {./broy3d_n1000_nops.dat};
	    \addplot+[line width= 0.3mm,mark=none, black] table {./broy3d_n1000_ps.dat};    
	\legend{unstructured, structured (ps)};
	\end{semilogyaxis}
	\end{tikzpicture}
	\caption{\label{br1000-f}The evolution of $f$ as a function of the number of (complete)
  evaluations for {\tt BROYDEN3D} for $n=1000$}
	\end{figure}

The rate of decrease of the structured BFO is very clearly much faster than
that of the unstructured variant. Moreover, the two passes of the
structure-exploiting algorithm are very noticeable, the first pass being
significantly faster and yielding most of the final decrease.  We observe that
the length of the second pass increases with size, which is expected because
the space to sample in the second pass increases in dimension.  However,
because the first pass already made such good progress, the effort spent in
the second pass remains acceptable, at least for the problem sizes considered
here. By contrast, the performance of the standard unstructured version of BFO
quickly degrades with size: for small dimensions, one notices the
staircase-like decrease which is typical of pattern search methods.  In this
test problem, the first poll step is unsuccessful, leading to a full $2n$
function evaluations before the stepsize is reduced.  While this can be
acceptable for small $n$, it becomes more problematic as $n$ grows. For larger
problems, the structured variant has already terminated before the first poll
step is completed in the unstructured method. For the instance in 10
variables, one also notices (in the left picture of Figure~\ref{br10-50-f})
that it can be beneficial to reapply the structured poll-step mechanism of the
first pass after a sufficient decrease in the second pass (cfr. Figure~\ref{fig:example}).

The behaviour of the structured variants on problem {\tt BDEXP} finally merits
a comment.  This problem features a very flat objective function near the
solution and, while the objective function value is decreased quickly, the
mechanism of the pattern search method takes many iterations to declare
optimality and terminate.

\subsection{Exploiting structured models in the search step}

We now turn to illustrating the performance which can be obtained using the
BFOSS model library described above, and also compare it with the pure
(structured and unstructured) sampling strategy of BFO alone.  In the experiments
reported next, we use fully quadratic models, the Mor\'{e}-Sorensen
method for maximizing the Lagrange polynomials in the trust-region and the 
truncated conjugate-gradient algorithm to solve \req{trbox}.

In this section we use performance and data profiles 
  \cite{DolaMore02,MoreWild09} to compare different variants of BFO.
For this purpose, we measure performance in terms of numbers of
  full objective function evaluations necessary for termination that is
  declared when the following condition holds
\begin{equation}\label{convtest}
f(x_0)-f(x) \ge (1-\tau)(f(x_0) - f_*).
\end{equation}
Here $x_0$ is the starting point for the problem, $x$ is the solution
returned by a solver, $f_*$ is computed for each problem as the smallest value
of $f$ obtained by any solver within a given number $\mu_f$ ($\mu_f = 100000$ in our tests) of function
evaluations, $\tau \in[0,1]$ is a tolerance that represents the the percentage
decrease from the starting value $f(x_0)$ (we used the standard value $\tau = 10^{-4}$). The stopping criterion \req{convtest} is
suggested in \cite{MoreWild09} to generate profiles and differs from the default 
stopping criterion on the minimum step size used in BFO by default.

In Figure~\ref{perf_small_f} four variants of BFO are compared on the small test. 
These variants are
\begin{description}
\item[unstructured:]   the standard unstructured BFO algorithm without using the
     BFOSS models,
\item[models:] the standard unstructured BFO algorithm where a full-dimensional
     BFOSS search step is attempted at every iteration,
\item[ps:]    the version of BFO using the coordinate partially-separable structure,
     but without using the BFOSS models,
   \item[ps \& models:] the version of BFO using the coordinate partially-separable
     structure, using BFOSS models for each element function at every iteration.
\end{description}

\comment{
\begin{figure}[htbp]
\centerline{
\includegraphics[height=5cm,width=7cm]{small_all.eps}
}
\caption{\label{perf_small_f}Performance profile for unstructured/structured variants,
  with or without models (small test set)}
\end{figure}
}

\begin{figure} \centering
		\begin{tikzpicture}
		\begin{axis}[width=.5 * \linewidth, 
		title = {Performance profiles with $\tau=10^{-4}$},
		xlabel = {$t$}, 
		height = .25 \textheight,
		legend style={at={(0.5,-0.41)},anchor=north},
		legend cell align={left},
		xmax = 10,
		xmin = 1,
		ymax = 1.1,
		ymin = 0
		]
\addplot+[thick, mark=none, magenta] table {./pp_small_uns_0.0001.dat};		
\addplot+[thick, mark=none, red, dashed] table {./pp_small_unsmod_0.0001.dat};
\addplot+[thick, mark=none, black] table {./pp_small_ps_0.0001.dat};		
\addplot+[thick, mark=none, blue, dashed] table {./pp_small_psmod_0.0001.dat};

\legend{unstructured,  models, ps, ps \& models};
\end{axis}
\end{tikzpicture}
		\begin{tikzpicture}
		\begin{axis}[width=.5 * \linewidth, 
		title = {Data profiles with $\tau=10^{-4}$},
		xlabel = {$\nu$}, 
		height = .25 \textheight,
		legend style={at={(0.5,-0.41)},anchor=north},
		legend cell align={left},
		xmax = 1500,
		xmin = 0,
		xtick = {0,400,...,2000},
		ymax = 1.1,
		ymin = 0
		]
\addplot+[thick, mark=none, magenta] table {./dp_small_uns_0.0001.dat};		
\addplot+[thick, mark=none, red, dashed] table {./dp_small_unsmod_0.0001.dat};
\addplot+[thick, mark=none, black] table {./dp_small_ps_0.0001.dat};		
\addplot+[thick, mark=none, blue, dashed] table {./dp_small_psmod_0.0001.dat};

\legend{unstructured,  models, ps, ps \& models};
\end{axis}
\end{tikzpicture}
 \caption{\label{perf_small_f}Performance and data profiles for unstructured/structured variants,
  with or without models (small test set).}
\end{figure}
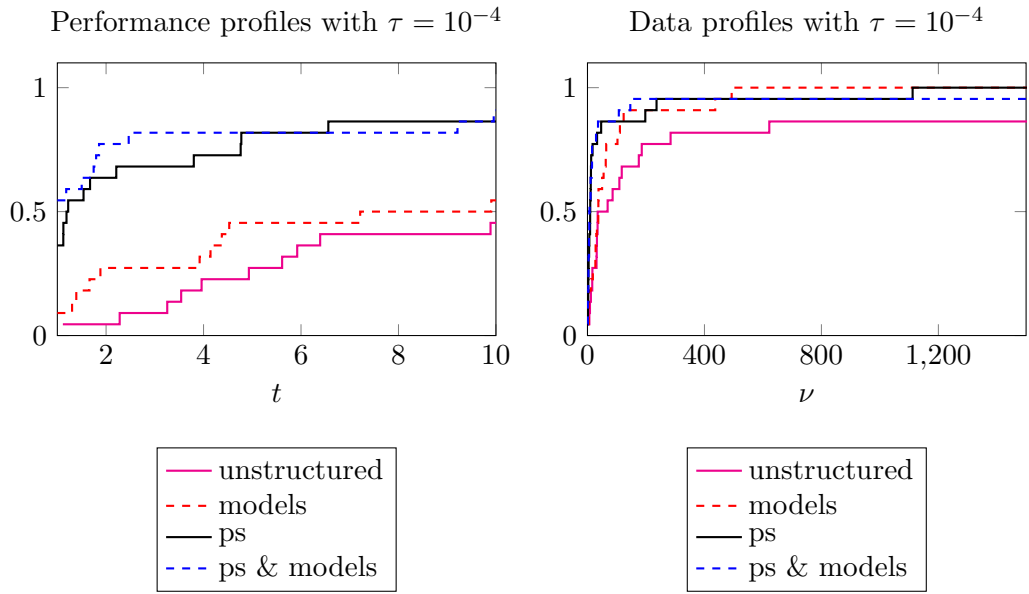

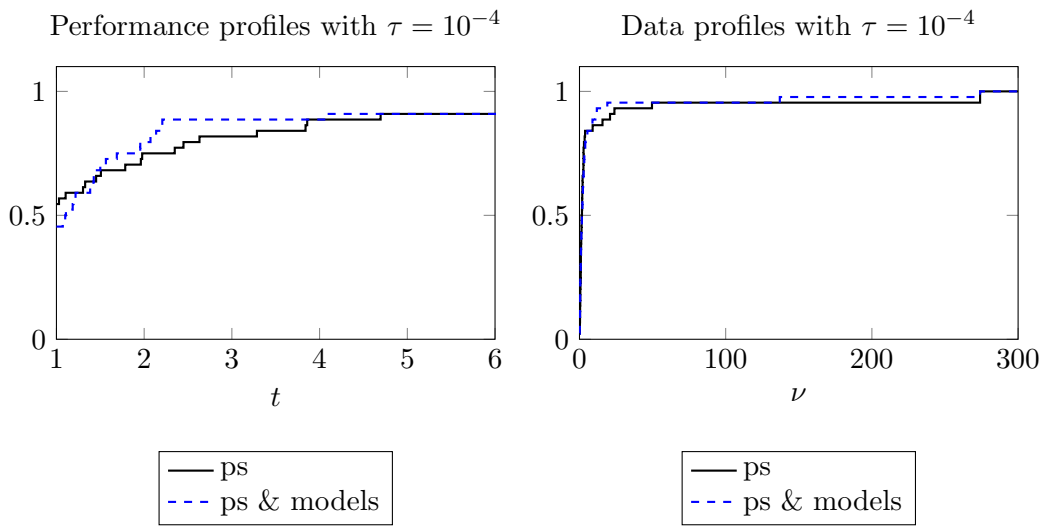
\begin{figure}\centering
		\begin{tikzpicture}
		\begin{axis}[width=.5 * \linewidth, 
		title = {Performance profiles with $\tau=10^{-4}$},
		xlabel = {$t$}, 
		height = .25 \textheight,
		legend style={at={(0.5,-0.41)},anchor=north},
		legend cell align={left},
		xmax = 6,
		xmin = 1,
		ymax = 1.1,
		ymin = 0
		]
\addplot+[thick, mark=none, black] table {./pp_smallish_ps_0.0001.dat};		
\addplot+[thick, mark=none, blue, dashed] table {./pp_smallish_psmod_0.0001.dat};

\legend{ps, ps \& models};
\end{axis}
\end{tikzpicture}
		\begin{tikzpicture}
		\begin{axis}[width=.5 * \linewidth, 
		title = {Data profiles with $\tau=10^{-4}$},
		xlabel = {$\nu$}, 
		height = .25 \textheight,
		legend style={at={(0.5,-0.41)},anchor=north},
		legend cell align={left},
		xmax = 300,
		xmin = 0,
		xtick = {0,100,...,300},
		ymax = 1.1,
		ymin = 0
		]
\addplot+[thick, mark=none, black] table {./dp_smallish_ps_0.0001.dat};		
\addplot+[thick, mark=none, blue, dashed] table {./dp_smallish_psmod_0.0001.dat};

\legend{ps, ps \& models};
\end{axis}
\end{tikzpicture}
 \caption{\label{perf_smallish_f}Performance and data profiles for structured variants,
  with or without models (smallish test set).}
\end{figure}

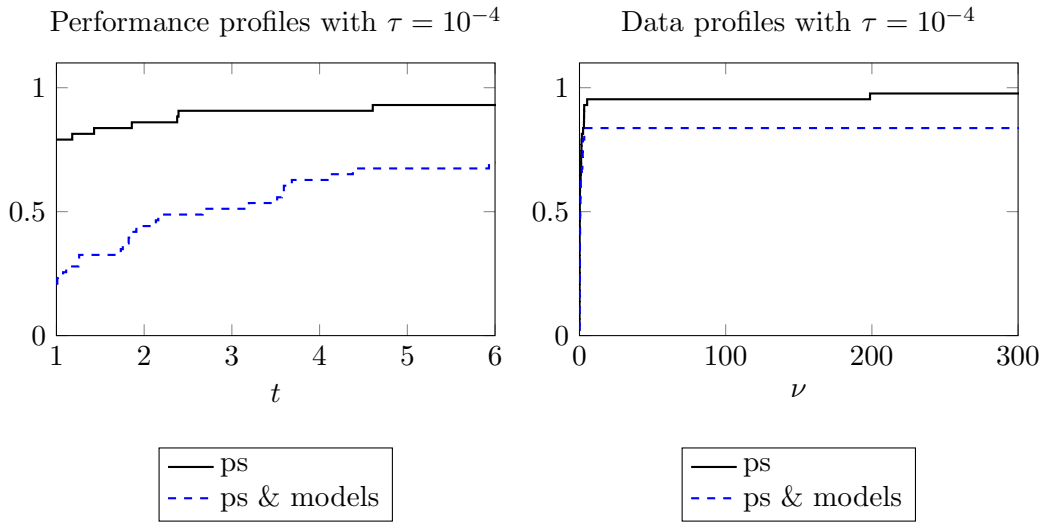
\begin{figure}\centering
		\begin{tikzpicture}
		\begin{axis}[width=.5 * \linewidth, 
		title = {Performance profiles with $\tau=10^{-4}$},
		xlabel = {$t$}, 
		height = .25 \textheight,
		legend style={at={(0.5,-0.41)},anchor=north},
		legend cell align={left},
		xmax = 6,
		xmin = 1,
		ymax = 1.1,
		ymin = 0
		]
\addplot+[thick, mark=none, black] table {./pp_medium_ps_0.0001.dat};		
\addplot+[thick, mark=none, blue, dashed] table {./pp_medium_psmod_0.0001.dat};

\legend{ps, ps \& models};
\end{axis}
\end{tikzpicture}
		\begin{tikzpicture}
		\begin{axis}[width=.5 * \linewidth, 
		title = {Data profiles with $\tau=10^{-4}$},
		xlabel = {$\nu$}, 
		height = .25 \textheight,
		legend style={at={(0.5,-0.41)},anchor=north},
		legend cell align={left},
		xmax = 300,
		xmin = 0,
		xtick = {0,100,...,300},
		ymax = 1.1,
		ymin = 0
		]
\addplot+[thick, mark=none, black] table {./dp_medium_ps_0.0001.dat};		
\addplot+[thick, mark=none, blue, dashed] table {./dp_medium_psmod_0.0001.dat};

\legend{ps, ps \& models};
\end{axis}
\end{tikzpicture}
 \caption{\label{perf_medium_f}Performance and data profiles for structured variants,
  with or without models (medium test set).}
\end{figure}

This figure shows that, for small problems, the combined use of models and structure
is the best algorithmic choice, but also indicates that using structure without models
is clearly preferable to using models without structure.  An interpretation of this
observation is that providing global information on the problem (structure) outperforms
approximating local one (models). The situation is less clear when the size of the problems
increases, as is shown in Figures~\ref{perf_smallish_f} and \ref{perf_medium_f}, where one compares the
performance of the 'ps' and 'ps \& models' variants on the smallish
and medium test sets (the two other variants fail for a large
proportion of the smallish problems). Profiles on the large test-set have not been 
generated since 'ps' was the only variant that could solve the test set in
the maximum time allowed (8 hours).

\comment{
\begin{figure}[htbp]
\centerline{
\includegraphics[height=5cm,width=7cm]{small&smallish_ps_psm.eps}
\includegraphics[height=5cm,width=7cm]{medium_ps_psm.eps}
}
\caption{\label{perf_larger_f}Performance profile for structured variants, with
  and without models (left: small and smallish test sets, right: medium test set)}
\end{figure}
}

In Figures~\ref{perf_smallish_f} and \ref{perf_medium_f}, we see that the relative
advantage obtained by the use of models for small
problems progressively vanishes to disappear completely when the problem size grows.

While comparing the number of function evaluations, as we have done above, is
most natural for derivative-free problems (the cost of an evaluation in real
world applications often dominating that of all computations internal to the
algorithm), it is also interesting to consider memory usage and internal
computing effort. The pure unstructured variant is clearly the most economical
from both points of view (but at the price of being the less efficient in
function evaluations). Because of the need to store the model itself, memory
usage can significantly increase for the variant using unstructured models (at
least for fully quadratic ones as discussed above).  This is also the case for
the structured variants because they have to store, analyze (once) and exploit
the structure, which requires using additional pointers and lists.  However,
the additional memory necessary to use models in the structured case remains
typically modest, as only a collection of small matrices needs being stored.
We finally note that variants using models require a significantly higher
internal computational effort, mostly in the solution of the many trust-region
subproblems involved in managing the interpolation set and computing the
search step.
This is apparent in Figure~\ref{perf_medium_f}
and in Table~\ref{results-full-t} in Appendix \ref{app:res} where
we report the number of
function evaluations taken by all variants on all test sets  using the
internal default stopping criterion. 
It is clear from the table 
that the worse reliability of the 'ps \&
models' variant is nearly entirely due to exhausting the maximum allowed
cpu-time (symbol '$\dagger$' in the table).
 Numerical experiments were also carried out using linear and
diagonal models in the search step (both in the structured and unstructured
cases). Unfortunately, it was observed that, despite the solution of the involved trust-region problem being considerably cheaper than 
when using quadratic models, the use of the simpler models resulted in a
considerable increase of the number of function evaluations, an thus worse overall performance.

\numsection{Conclusions and perspectives}\label{concl_s}

We have introduced a structured version of derivative-free random pattern
search algorithms which is able to exploit coordinate partially separable
structure (typically associated with sparsity) present in unconstrained and
bound-constrained optimization problems.  This techniques improves performance
by orders of magnitude and makes it possible to solve problems that otherwise
are totally intractable by standard derivative-free methods.

We have also described a library of interpolation-based modelling tools which
can be associated to the structured or unstructured versions of the initial 
pattern search algorithm. 
For problems of small or moderate size, the use of the library further enhances performance,
especially when associated with structure. For larger problems the internal computing costs
increase and, while still reducing the number of function evaluations, the use of the
library may require a sometimes unrealistic computing time, in particular if the problem
is unstructured.

In comparing the benefits of using problem structure in a poll step and
building local models using interpolation techniques, we have concluded that the former is likely to
be more efficient, in particular for larger problems, even for
structure-exploiting models.

A new release of the Matlab BFO package\footnote{Although this is unrelated to
  the subject of the present paper, we emphasize that the new BFO release also
  features the important capability of handling categorical
  variables, as well as the new training strategies discussed in
  \cite{PorcToin17c}, see Appendix \ref{app:BFO2}.} 
  featuring both use of structure and modelling
tools (as discussed in this paper) is now freely available online from \\

\centerline{\fbox{\tt https://github.com/m01marpor/BFO}}
\vskip 10pt
\noindent 
and the main new features are briefly sketched in Appendix \ref{app:BFO2}.

The selective
use of structured models in conjunction with structured poll steps remains an
attractive option. However, the criteria defining the circumstances in which
interpolation models should be used, if at all, need further investigation.
Many other topics also merit research, including the design of ``grayer''
optimization tools which could exploit derivatives available for part of the
problem while adapting the techniques described here for the rest.

{\footnotesize

\section*{Acknowledgements}

Both authors are indebted to the University of Florence for its support while the present research was being conducted.
   
The first author is member of the {\em Gruppo Nazionale per il Calcolo Scientifico} (GNCS) of the Istituto Nazionale di Alta Matematica (INdAM) and this work was partially supported by INdAM-GNCS under Progetti di Ricerca 2019 and 2020.

}

\clearpage

\appendix
\newpage

\section{Full results}\label{app:res}
Table \ref{results-full-t} shows  the number of objective function evaluations required by
  the structured/unstructured versions of BFO with and without interpolation model
  for the problem instances of  Table~\ref{probs-t} ($\infty$ meaning that more
  than 100000 evaluations were needed, $\dagger$ meaning that cpu time exceed
  the maximum time, 8 hours, allowed). Moreover, a blank entry
    means that the solution of the corresponding problem instance was not
    attempted either because  the solver version failed to solve a smaller
    version of the same problem (symbols $\dagger$ and $\infty$) or because
    the problem size was considered too large for the considered solver. In
    particular, the unstructured/model implementation was used for solving small problems only and the partially separable/model implementation was not applied to solve large problems).

\begin{table}
\begin{center}
\tiny
\begin{tabular}{l|r|rrrr|l|r|rrrr}
               &       &\multicolumn{2}{c}{unstructured} & \multicolumn{2}{c}{partially separable} &
               &       &\multicolumn{2}{c}{unstructured} & \multicolumn{2}{c}{partially separable} \\
Pb             &  $n$  & no mod  & model  & no mod &  model     &  Pb     & $n$   & no mod  &  model & no mod &  model     \\
\hline
               &    10 &     962   &  365   &     79   &   148    &               &    10 &    1199   &  1227  &    101   &   153  \\
               &    50 &   11859   &        &     91   &   122    &               &    50 &  $\infty$ &        &     70   &   381    \\
               &   100 &   36085   &        &     97   &   214    &               &   100 &           &        &     71   &   345   \\
{\tt ARWHEAD}  &   500 & $\infty$  &        &    146   &   659    & {\tt BDARWHD} &   500 &           &        &     81   &   691   \\
               &  1000 &           &        &    194   &  1015    &               &  1000 &           &        &     76   &   869   \\
               &  5000 &           &        &    389   &          &               &  5000 &           &        &     77   &          \\
               & 10000 &           &        &    618   &          &               & 10000 &           &        &     77   &          \\
\hline

               &    10 &  21122    &  222   &   1661   &   375    &               &    10 &    2468   &  883   &    290   &    345         \\
               &    50 & $\infty$  &        &  13218   &   828    &               &    50 &  $\infty$ &        &    301   &   1189      \\
               &   100 &           &        &  34408   &   156    &               &   100 &           &        &    298   &    457           \\
{\tt BDEXP}    &   500 &           &        &$\infty$  &   137    & {\tt BDQRTIC} &   500 &           &        &    393   &    373         \\
               &   1000&           &        &          &   132    &               &  1000 &           &        &    542   &    503       \\
               &   5000&           &        &          &          &               &  5000 &           &        &   1150   &             \\
               &  10000&           &        &          &          &               & 10000 &           &        &   1480   &            \\
\hline
               &    10 & $\infty$  &$\infty$&    275   &    70    &               &    10 &    1225   &  456   &    308   &    96          \\
               &    50 &           &        &    275   &    85    &               &    50 &   22244   &        &    199   &   118              \\
               &   100 &           &        &    275   &   123    &               &   100 &   74758   &        &    273   &   151         \\
{\tt BEALES}   &   500 &           &        &    275   &   203    &{\tt BROYDN3D} &   500 & $\infty$  &        &    304   &   215         \\
               &  1000 &           &        &    275   &   234    &               &  1000 &           &        &    370   &   303                  \\
               &  5000 &           &        &    300   &          &               &  5000 &           &        &    640   &                        \\
               & 10000 &           &        &    325   &          &               & 10000 &           &        &    675   &                        \\
\hline 
               &    10 &    6171   & 1224   &  15773   &   567     &              &    16 &     268   &   126  &    265   &  107           \\
               &    50 & $\infty$  &        &   4325   &   816     &              &    64 &   30452   &        &    221   &  248                \\
               &   100 &           &        &   2788   &  1798     &              &   144 &  $\infty$ &        &    442   &  494            \\
{\tt BROWNAL6} &   500 &           &        &   2682   & $\dagger$ &{\tt CONTACT} &   400 &           &        &   1009   &  1295                  \\
               &  1000 &           &        &   2847   & $\dagger$ &              &   900 &           &        &   1761   &  $\dagger$                    \\
               &  5000 &           &        &   3059   &           &              &  2500 &           &        &   3589   &                        \\
               & 10000 &           &        &   3118   &           &              &  4900 & $\infty$  &        &   5952   &                       \\
    
\hline
               &    15 &    2358   &  558   &    185   &   113    &               &    15 & 17702     &  1405  &    185   &   152          \\
               &    51 &   17357   &        &    189   &   101    &               &    51 & $\infty$  &        &    186   &   195               \\
               &   102 &   55748   &        &    240   &    96    &               &   102 &           &        &    184   &   284              \\
{\tt DIXMAANA} &   501 &  $\infty$ &        &    226   &   313    &{\tt DIXMAANI} &   501 &           &        &    230   &   1231               \\
               &  1002 &           &        &    375   &   869    &               &  1002 &           &        &    265   &   1487              \\
               &  5001 &           &        &    310   &          &               &  5001 &           &        &    709   &                         \\
               & 10002 &           &        &    756   &          &               & 10002 &           &        &    798   &                         \\
\hline
               &    10 &    1567   &  629   &    143   &   83     &               &    10 &   83101   &  4570  &    257   &   2297         \\
               &    50 &   30400   &        &    155   &   94     &               &    50 & $\infty$  &        &    181   &   2183                  \\
               &   100 &  $\infty$ &        &    157   &  129     &               &   100 &           &        &    191   &   2712       \\
{\tt ENGVAL}   &   500 &           &        &    159   &  158     &{\tt FREUROTH} &   500 &           &        &    185   &   $\dagger$                    \\
               &  1000 &           &        &    159   &  218     &               &  1000 &           &        &    192   &   $\dagger$                 \\
               &  5000 &           &        &    159   &          &               &  5000 &           &        &    233   &                        \\
               & 10000 &           &        &    159   &          &               & 10000 &           &        &    317   &                        \\
\hline
               &    11 &    9931   & 2094   &    131   &  101     &               &    24 &     734   &    65  &    101   &    19           \\
               &    51 & $\infty$  &        &    152   &  113     &               &    84 &   27567   &        &    301   &   779                 \\
               &   101 &           &        &    167   &  108     &               &   180 &  $\infty$ &        &    427   &   37               \\
{\tt HELIX}    &   501 &           &        &    298   &  122     &{\tt JNLBRNG1} &   544 &           &        &    944   &   86               \\
               &  1001 &           &        &    388   &  152     &               &  1012 &           &        &   1393   &   92             \\
               &  5001 &           &        &    883   &          &               &  1984 &           &        &   1909   &                         \\
               & 10001 &           &        &   1070   &          &               &  5304 &           &        &   1799   &                         \\
\hline
               &    16 &     318   &  117   &    306   &   76     &               &    12 &    7154   &   391  &   7010   &  117       \\
               &    64 &   23692   &        &    461   &  243     &               &    52 &   22732   &        &     60   &  242               \\
               &   144 &  $\infty$ &        &   1065   &  516     &               &   102 &    1923   &        &     47   &  263         \\
{\tt LMINSURF} &   400 &           &        &   3301   & 1580     &{\tt MOREBV}   &   502 &    9001   &        &     47   &   55                   \\
               &   900 &           &        &   8506   &$\dagger$ &               &  1002 &   18001   &        &     47   &   56                    \\
               &  2500 &           &        &  29961   &          &               &  5002 &   90001   &        &     47   &                         \\
               &  4900 &           &        &  $\infty$&          &               & 10002 &   $\infty$&        &     47   &                         \\
\hline
               &    16 &    415    &   94   &   412    &   76     &               &    13 &    1480   &  1397  &    177   &   297             \\
               &    64 &  29409    &        &   579    &  358     &               &    39 &   61772   &        &    225   &   347             \\
               &   144 &  $\infty$ &        &  1082    &  623     &               &   130 &  $\infty$ &        &    625   &   523                 \\
{\tt NLMINSRF} &   400 &           &        &  8853    & 1824     &{\tt NZF1}     &   650 &           &        &    684   &  1476                  \\
               &   900 &           &        &  3412    & $\dagger$&               &  1300 &           &        &    667   &  2185               \\
               &  2500 &           &        &  31619   &          &               &  6500 &           &        &    946   &                         \\
               &  4900 &           &        & $\infty$ &          &               & 13000 &           &        &   2228   &                         \\
\hline
               &    20 &   20605   &  2110  &    716   &  776     &               &    10 &   13241   & 12229  &    361   &   299      \\
               &    52 & $\infty$  &        &    824   &   81     &               &    50 & $\infty$  &        &    361   &   534                \\
               &   100 &           &        &    849   &   84     &               &   100 &           &        &    384   &   729          \\
{\tt POWSING}  &   500 &           &        &    988   &   94     &{\tt ROSENBR}  &   500 &           &        &    436   &  1531      \\
               &  1000 &           &        &   1036   &   98     &               &  1000 &           &        &    461   &  2334  \\
               &  5000 &           &        &   1164   &          &               &  5000 &           &        &    636   &             \\
               & 10000 &           &        &   1148   &          &               & 10000 &           &        &    736   &              \\
\hline
               &    10 &    3073   &   139  &    440   &   17     &               &    20 & $\infty$  & 32414  &   1609   &   1538     \\
               &    50 & $\infty$  &        &    316   &   20     &               &    40 &           &        &   1747   &   1985               \\
               &   100 &           &        &    314   &   23     &               &   200 &           &        &   1924   &   4210          \\
{\tt TRIDIA}   &   500 &           &        &    345   &   30     &{\tt WOODS}    &   500 &           &        &   2241   &   $\dagger$        \\
               &  1000 &           &        &    293   &   33     &               &  1000 &           &        &   2927   &   $\dagger$      \\
               &  5000 &           &        &    277   &          &               &  5000 &           &        &   4515   &\\
               & 10000 &           &        &    278   &          &               & 10000 &           &        &   5002   &            \\
\hline
\end{tabular}
\caption{\label{results-full-t}The number of objective function evaluations required by
  the structured/unstructured versions of BFO with and without interpolation model
  for the problem instances of  Table~\ref{probs-t}. }
\end{center}
\end{table}

\section{Details for problem {\tt CONTACT}}\label{app:contact}

The problem {\tt CONTACT} is a bilinear finite-element discretization of a
membrane contact problem defined on the unit square of $\Re^2$.  This square
is discretized in $q = (\sqrt{n}-1)^2$ element squares and the surface of the
membrane ``above'' this element square is given by
\[
\frac{1}{q}\sqrt{1+ (x_{SW}-x_{NE})^2 + (x_{SO}-x_{NW})^2},
\]
where $x_{SW}$, $x_{NE}$, $x_{SO}$ and $x_{NW}$ give
the height of the membrane at the four corners of the element square.  The membrane
is fixed on the boundaries of the unit square to the value of the function
\beqn{c-boundary}
b(x,y) = 1 + 8x + 4y + 3 \sin( 2\pi x )\sin( 2 \pi y ).
\eeqn
In addition, the membrane has to lie above a square flat obstacle of height 10
positioned at $[0.4,\,0.6]^2$.  The starting point is the projection of the
vector given by the values of \req{c-boundary} at each element square's corner
onto the feasible domain.

\section{Details for problem {\tt NZF1}}\label{app:nzf1}

The {\tt NZF1} problem is a variable dimension version of the eponymous
nonlinear least-squares problem defined in \cite{PricToin06} as an example
illustrating the concept of partially separable functions.  If
$n = 13\ell$, it is defined by
\[
\begin{array}{lcl}
f(x_1, \ldots, x_{n})
& = & \bigsum_{i=1}^\ell \Bigg[ \Big( 3x_i - 60 + \sfrac{1}{10}( x_{i+1}-x_{i+2} )^2 \Big)^2 \\
&   & \ms\ms  + \left( x_{i+1}^2 + x_{i+2}^2 + x_{i+3}^2(1+x_{i+3})^2 + x_{i+6}
                 + \bigfrac{x_{i+5}}{(1 + x_{i+4}^2 + \sin\left( \frac{x_{i+4}}{1000}\right)}  \right)^2\\
&   &\ms\ms     +\Big( x_{i+6} + x_{i+7} - x_{i+8}^2 + x_{i+10} \Big)^2 \\
&   &\ms\ms     +\Big( \log( 1 + x_{i+10}^2 ) + x_{i+11} - 5x_{i+12} + 20\Big)^2 \\
&   &\ms\ms     +\Big( x_{i+4} + x_{i+5} + x_{i+5}*x_{i+9} + 10x_{i+9} - 50 \Big)^2\Bigg]\\
&   & + \bigsum_{i=1}^{\ell-1} (x_{i+6}-x_{i+19})^2.
\end{array}
\]
The starting point is the vector of all ones.

\section{BFO 2.0 and its new features}\label{app:BFO2}

Release~2.0 of the Matlab BFO package is a major upgrade from Release~1 and
includes several important new problem-oriented possibilities. Beyond the
BFOSS library of model-based search steps described in Section~\ref{model-s}
and the exploitation of CPS structure described in Section~\ref{pattern-s},
BFO~2.0 also supports the following new problem  features.
\begin{description}
\comment{
\item[Using coordinate partially-separable problem structure.]
  A large number of (often relatively large) optimization problems have some
  underlying structure, and BFO can be made to exploit this structure to
  considerable advantage.
 
  The first case is when the objective function can be expressed as a sum
  of the form
  \[
  f(x) = \sum_{i=1}^p f_i(x)
  \]
  in which case we say that the objective function is in \emph{sum form}.  The
  functions $f_i(x)$ are the called \emph{element functions}.  A first advantage of
  problems in sum form lies in the possibilities 
  for the user to define a structure exploiting search step function (for
  instance by building individual models for each of the element functions).
  Indeed, if the objective function is defined in sum form,
  BFO will automatically maintain an evaluation history for each individual
  element function, and will pass this information to the user when calling
  the user-defined search step (using BFOSS for instance).
   
  The main advantage of problems whose objective function is in sum form is
  the specialization of this structure to the very important case where the
  problem is "coordinate partially separable" (CPS) or "sparse". In this extremely
  frequent case (e.g. when the objective function is related to a problems with
  distinct interconnected blocks, or is resulting from the discretization of a
  continuous problem), the objective function has the form
\[
f(x) = \sum_{i=1}^p f_i(x_i),
\]
  where now $x_i$ only contains a subset of the problem's variables. A
  sum-form function is considered coordinate partially separable when 
  the maximal dimension of one of the $x_i$ is (often much) smaller than the
  total problem dimension. For example, the function
\[
 f( [ x(1) x(2), x(3) ]) = {\tt norm}( [ x(1) x(2) ], 2 )^2 + {\tt norm}( [ x(2) x(3) ], 'inf' )
\]
   is coordinate partially separable with 2 elements of domains defined by the
   two vectors of variable indices $[ 1 2 ]$ and $[ 2 3 ]$. Coordinate
   partially-separable functions as defined above are often called "sparse"
   because their Hessians (when they exist) are sparse matrices whose maximal
   dense principal submatrices are defined by the indices defining the $x_i$.

   The use of any underlying coordinate partially separable structure results in
   very substantial gains in the number of evaluations of $f(x)$.
   Even if the use of this feature may not be critical for small problems,
   its use for problems of moderate or large size is often essential,
   especially if the cost of evaluating the element functions $f_i(x)$ is high.
   The gains in efficiency (evaluation counts) may typically be of several orders
   of magnitude (the amount of computation and storage internal to the
   algorithm therefore increases relatively to the number of evaluations).
   (As demonstrated in Section~\ref{numerics_s}, the evaluation gains can be combined 
   with those resulting from a intelligent search-step step
   strategy due to the sum-form of the problem.
\item[The BFOSS library of model-based search steps.]
}
   
\item[Categorical variables.]
  In addition to standard continuous and discrete variables, BFO now supports
  the use of categorical variables. Categorical variables are 
  unconstrained non-numeric variables whose possible 'states' are defined
  by strings (such as 'blue'). These states are not implicitly ordered,
  as would be the case for integer or continuous variables. As a consequence,
  the notion of neighbourhood of a categorical variable is entirely
  application-dependent, and has to be supplied, in one form of two possible forms,
  by the user. Moreover, the 'vector of variables' is no longer a standard
  numerical vector when categorical variables are present, but is itself a
  \emph{vector state} defined by a value cell array of size $n$ (the problem's
  dimension), whose $i$-th  component is either a  number
  when variable $i$ is not categorical, or a string defining the current state
  of the $i$-th (categorical) variable.  For example, such a 4-dimensional vector
  state can be given by the value cell array
\begin{quote}
      \{\{ 'blue', 3.1416, 'magenta', 2 \}\}.
\end{quote}
  Variable $i$ is declared to be categorical by specifying xtype($i$) = 's'. If
  a problem contains at least one categorical variable, it is called a
  categorical problem and optimization is carried on vector states (instead of
  vectors of numerical variables).  As a consequence, the starting point
  and the returned best minimizing point are vector states, and the objective
  function's value is computed at vector states (meaning that the argument of
  the function $f$ is a vector state). 
  
  The user must specify the application-dependent neighbours (also
  called categorical neighbourhoods) of each given vector state with respect
  to its categorical variables. This can be done in two mutually exclusive
  ways.
  \begin{enumerate}
  \item The first is to specify \emph{static neighbourhoods}. This is done by
      specifying, for each categorical variable, the complete list of its
      possible states.  The
      neighbourhood of a given vector state $vs$ wrt to categorical variable $j$
      (the hinge variable) then consists of all vector states that differ
      from $vs$ only in the state of the $j$-th variable, which takes all possible
      values different from $vs(j)$. In this case, all variables of the problem
      retain their (initial) types and lower and upper bounds.
   \item The second is to specify \emph{dynamical neighbourhoods}.  This more flexible
      technique is used by specifying a user-supplied function whose purpose 
      is to compute the neighbours of the vector state $vs$ 'on the fly', when
      needed by BFO. At variance with the static neighbouhood case, the
      variable number and types, as well as lower and upper bounds of the
      neighbouring vector states (collectively called the 'context') may be
      redefined within a framework  defined by a few simple rules.

      In effect, this amounts to specifying the neighbouring nodes in a (possibly
      directed) graph whose nodes are identified by the list, types, bounds and
      values of the variables. As a consequence, the user-supplied definition of
      the neighbour(s) of one such node may need to take the values of all variables
      into account. Of course, for the problem to make sense, it is still required
      that the objective function can be computed for the new neighbouring vector
      states and that its value is meaningfully comparable to that at $vs$.
  \end{enumerate}
  The very substantial flexibility allowed by this mechanism of course comes
  at the price of the user's full responsibility for overall coherence.

\item[Performance and data profile training strategies.]
  Because BFO is a \emph{trainable package} (meaning that its internal
  algorithmic constants can be trained/optimized by the user to optimize its performance
  on a specific problem class), it needs to define training strategies which
  allow to decide if a particular option is better than another.  Release~1 of
  BFO included the natural ``average'' training criterion (quality is measured
  by the average number of function evaluations on the class) and a robust
  variant of the same idea (see \cite{PorcToin17} for details).  Release~2.0 now includes two new training
  strategies (introduced in \cite{PorcToin17c}):
  \begin{description}
  \item[Performance profiling.]
    When this training option is selected, the performance of two algorithmic variants
    (i.e.\ versions of BFO differing by the value of their internal
    algorithmic parameters) are compared using the well-known performance
    profile methodology \cite{DolaMore02,DolaMoreMuns06}.
  \item[Data profiling.]
    This option is similar to performance profiling, but uses data profiles
    \cite{MoreWild09} instead of performance profiles to compare two variants.
  \end{description}
  These new options correspond more closely to the manner in which
  derivative-free packages are compared in the optimization literature.
\end{description}


\end{document}